\documentclass{article}

\usepackage[pdftex]{graphicx}

\usepackage{blkarray}

\usepackage{TeX_aux_file_1_Macros_Colas_Oct2010}

\newtheorem*{lemmaSvdP}{Lemma S-vdP}
\newtheorem*{PropositionA}{Proposition A}
\newtheorem*{PropositionB}{Proposition B}
\newtheorem*{theoremintro}{Theorem} 
 
\newtheorem*{lemmaintro}{Lemma}

\title{
{On the geometrization of a lemma of differential Galois theory}}

\author{%
Colas Bardavid\footnote{colas.bardavid@gmail.com, bardavid@imsc.res.in} \\
Institute of Mathematical Sciences\\
CIT Campus, Taramani \\
Chennai 600113 --- Tamil Nadu,
India
}

\begin{document}

\maketitle
 \begin{center}  
\rule{10 cm}{.5pt} \end{center}
\begin{small}
\begin{center}
\begin{minipage}{10cm}
\textbf{Abstract --} 
In this paper, we give a geometrization and a generalization of a lemma of differential Galois theory. This geometrization, in addition of giving a nice insight on this result, offers us the occasion to investigate several points of differential algebra and differential algebraic geometry. We study the class of simple \D-schemes and prove that they all have a coarse space of leaves.  Furthermore, instead of consi\-dering schemes endowed with one vector field, we consider the case of arbitrarilly large families of vector fields. This leads us to some develop\-ments in differential algebra, in particular to prove the existence of the trajectory in this setting but also to study  simple \D-rings.
\end{minipage}
\end{center}
\end{small}
 \begin{center}  
\rule{10 cm}{.5pt} \end{center}

\begin{center}
 \begin{minipage}{10cm}
\emph{2000 Mathematics Subject Classification}: 12H05, 13B05, 13N10
\smallskip \\
\emph{Keywords}: differential Galois theory, differential algebra,  simple rings,  trajectory, coarse space of leaves
\end{minipage}
\end{center}

\newpage
\tableofcontents

\newpage

\section{Introduction}

In their book \Cite{singvdp} on differential Galois theory, Singer and van der Put follow two approaches to define and study Picard-Vessiot extensions. For $k$  a differential field with an algebraically closed constant field of characteristic $0$, and $A\in M_n(k)$, regarding the linear differential equation
\begin{equation}
Y'=AY \tag{$*$}
\end{equation}
they call $R/k$  a \emph{Picard-Vessiot ring for \emph{($*$)}} if 
\begin{itemize}  \setlength{\itemsep}{-3pt} 
\item$R$ is a simple differential ring;
\item there exists $U\in GL_n(R)$ such that $U'=AU$;
\item $R=k[\,\text{entries of }U, \frac{1}{\det U}\,]$.
\end{itemize}
This definition is natural because fields have two analogues in the differential world: differential fields and simple differential rings. Recall that these latter are defined to be the differential rings $R$  whose only  differential ideals are $0$ and $R$. Then, they prove (Proposition 1.21 of \Cite{singvdp}) that $K/k$ is a Picard-Vessiot field for ($*$), is the usual sense, if and only if there exists a Picard-Vessiot ring $R$ such that $K=\Frac R$. The proof of this equivalence requires the following lemma (Lemma 1.23 of \Cite{singvdp}).

\begin{lemmaSvdP}
Let $k$ be a differential field, whose field of constants  $C$ is algebraically closed and of characteristic zero. One considers
the $k$-algebra $k[X_{ij}, \frac{1}{\det X_{ij}}]$ as a differential ring with the derivation defined by ${X_{ij}}'=0$. Then, the maps
$$I \longmapsto \left \langle I  \right \rangle\qquad \aand \qquad J \longmapsto J \cap C[X_{ij}, \frac
{1}{\det X_{ij}}]$$
are inverse bijections between the set of ideals $I$ of $C[X_{ij}, \frac{1}{\det X_{ij}}]$ and the set of  differentials ideals $J$ of $k[X_{ij}, \frac{1}{\det X_{ij}}]$.
\end{lemmaSvdP}

The main goal of this paper is to apply some tools developped in \Cite{Traj1} (namely leaves and trajectory) to give a geometrical insight on and to generalize this lemma.

\bigskip
Just as for commutative algebra with respect to algebraic geometry,  many (not to say all) results and definitions of differential  algebra gain generality and geometric clarity when expressed in a geometric framework. Such an approach to differential algebra, in the setting of schemes, can be found in \Cite{buium}, \Cite{gillet}, \Cite{dyckerhoff}, \Cite{Vojta}, \Cite{ProlongDiffAlg}. Note than in \Cite{EGA44}, Grothendieck writes 
\emph{``Nous passons sous
silence de nombreux d\'eveloppe\-ments, classiques en G\'eo\-m\'etrie diff\'e\-rentielle (connexions,
transformations infinit\'e\-si\-males associ\'ees \`a un champ de vecteurs, jets, etc.), bien que ces
notions s'\'ecri\-vent de fa\c{c}on particuli\`erement naturelle dans le cadre des sch\'emas. Nous
passons \'egalement sous silence ici les ph\'enom\`enes sp\'eciaux \`a la caract\'eristique \mbox{$p>0$}.''}\footnote{\emph{``We will not do here many developments that are classical for differential geometry
(connections,
infinitesimal transformations associated to a vector field, jets, etc.), although scheme theory is a very natural setting for these
concepts. Also, we will not touch on  pheno\-mena specific to characteristic $p>                                                                    0$.''}}

\bigskip

In this paper, we aim at contributing to this program.
When we first saw this statement, we were quite convinced that there were a geometric statement behind it and this was our starting point. Actually, this idea  has lead us to several developments on simple differential rings and schemes as well as on the existence and the computation of a coarse space of leaves of a scheme endowed with a vector field. Finally, our geometric and generalized version of Lemma S-vdP is the following (see Theorem \Ref{main_thm}):

\begin{theoremintro}
Let $\mathscr{S}$ be a simple \D-scheme  of characteristic zero and let  $\mathscr{X}$ be a $\mathscr{S}$-scheme without self-dynamics. Let $\mathscr{S} \longto C(\mathscr{S})$ be a coarse space of leaves of $\mathscr{S}$. Let $x$ be a $C(\mathscr{S})$-point of $\mathscr{X}$. Then, the maps
$$
\fonctionb{X_x}{{X}^{\chv{V}}}{y}{\mathrm{Traj}_\chv{V}(y)}
\qquad \aand \qquad
\fonctionb{X^\chv{V}}{X_x}{\eta}{\eta_{\mid x}}
$$ 
are inverse bijections between the fibre $X_x$ and the set $X^\chv{V}$ of leaves of $X$ under the action of $\chv{V}$.
\end{theoremintro}

\noindent By a \D-scheme, we mean a scheme endowed with a vector field (see Definition \Ref{DeltaScheme}).
Actually, this theorem says two things. First, if $\eta$ is any leaf (see Defi\-ni\-tion \Ref{Defi_leaf}) of $\mathscr{X}$, then $\eta$ intersects the fiber $X_x$ in a unique point $\eta_{\mid x}$. Note however that the dimension of $\eta_{\mid x}$ can be greater than $0$. Second, it says that 
$\Traj{V}(\eta_{\mid x})=\eta$ and ${\Traj{V}(p)}_{\mid x}=p$
for every leaf $\eta$ of $\mathscr{X}$ and every point $p$ of the fiber $X_x$.
 It is a generalization of  the lemma S-vdP along three directions:
\begin{itemize}
\item[---] 
 In Lemma S-vdP, the algebra $C[X_{ij}, \frac{1}{\det X_{ij}}]$ of the algebraic group $GL_{n,C}$ appears. We show that one can replace $GL_{n,C}$ by any scheme $X$.
\item[---] Second, we replace the base field $k$ of the lemma S-vdP by any scheme $S$ endowed with a \emph{simple} vector field. This means that the vector field $\chv{V}_S$  ``melts completely'' $S$ (see Definition \Ref{DScheme_simple}). This generalization is more difficult than the previous one. In Lemma S-vdP, the base is necessarily of  dimension zero while in our version the base can be arbitrarilly large. For instance, the line $\mathbf{A}^{1}_{\mathbf{C}}$ endowed with the usual vector field $\partial/\partial x$ is simple. Also, the scheme $C(\mathscr{S})$ (which plays the role of the constant field $C$) is not assumed to be algebraically closed.
 \item[---] Actually, the theorem we will prove is a little more general, since instead of dealing with schemes equipped with \emph{one} vector field, we will handle schemes endowed with any family of  vector fields. Herein, we are following the setting of Buium in \Cite{buium}. The usual generalization in this direction is to replace differential rings  by rings with $d$ commuting derivations ---so, here it is much more general.
\end{itemize}

\smallskip

For the reader to fully understand this result, it is needed to give some precisions:
\begin{itemize}
 \item[---] 
 First, as one can imagine, a coarse space of leaves $C(\mathscr{S})$ of some scheme $\mathscr{S}$ endowed with a vector field is a scheme that parametrizes in a ``cate\-gorical sense'' the leaves of $\mathscr{S}$ (see Definition \Ref{Def_coarse_space_of_leaves}). As we will see, it is a geometric analogue of the constant ring of a differential ring. Such coarses spaces do not exist in general, but we prove in Theorem \Ref{Existence_Coarse_Sp_Leaves_Qu_Simple} that  it always exists for (quasi-)simple schemes. Although the coarse space of leaves is not completely satisfactory, this tool here is sufficient. For a better approximation of the space of leaves, one option would be to follow the guideline of GIT \Cite{mumfordgroup} but we will not do it here. 
 \item[---] 
  Second, by a  $\mathscr{S}$-scheme without self-dynamics, we wean that  $\mathscr{X}$ has a model $X_{0}$ over $C(\mathscr{S})$. In the lemma S-vdP, the base $\mathscr{S}$ correponds to $k$, the coarse space of leaves $C(\mathscr{S})$ corresponds to $C$, the scheme $\mathscr{X}$ corresponds to $k[GL_{n}]$ and $X_{0}$ to $C[GL_{n}]$. In some way, saying that $\mathscr{X}/\mathscr{S}$ has no self-dynamics corresponds to the fact that one endows $k[X_{ij}, \frac{1}{\det X_{ij}}]$ with the derivation $X_{ij}'=0$.
\end{itemize}

\smallskip

This theorem reduces in the affine case to Proposition A (the very affine generalization of Lemma S-vdP) and Proposition B (another affine statement, which is needed to compare the model $X_0$ of $\mathscr{X}$ with the fiber $X_x$), whose proofs are given in Section {\ref{Affine_Statements}}. Proposition A reads as:
\begin{PropositionA}
 Let $K$ be a simple \D-ring and $C$ its field of constants. Let $A$ be a $C$-algebra, endowed with the zero derivation. Denote by $i : A \longto A\otimes_C K$ the canonical morphism. Then, the maps
$$I \longmapsto \left \langle i(I)  \right \rangle\qquad \aand \qquad J \longmapsto i^{-1}(J) $$
are inverse bijections between the set of ideals $I$ of $A$ and the set of  \D -ideals $J$ of $A\otimes_C K$.
\end{PropositionA}
\noindent  The proof of Proposition A follows partly the proof of Singer and van der Put but with two more important ideas. First, we proceed by induction on what we call the \emph{length} of an element of $A\otimes_C K$. The length of $x\in A\otimes_C K$ is the minimal length of an expression of $x$ as $x=\sum_i a_i\otimes \lambda_i$. The second idea is to consider, if $R$ is a \D-ring and if $x\in R$, the left ideal $ \text{EqDiffLin}_R(x)$ of $R[\partial]$ consisting in all $L$ such that $L\bullet x=0$. In a differential field $K$, for any nonzero element $x\in K$, there exists $\ell \in R$ such that $\ell\cdot x=1$. For simple \D-rings, one has a similar property, but less convenient to handle: for any nonzero element $x\in K$, there exists $L \in R[\partial]$ such that $L\bullet x=1$. 

\bigskip

This paper is organised as follows. Section 1 is this introduction. Section 2 is dedicated to introducing \D-schemes, leaves and trajectory. Whereas in \Cite{Traj1}, we were working with schemes endowed with one vector field, we construct here these tools for schemes endowed with any family of vector fields. In Section 3, we define simple and quasi-simple \D-schemes and prove that if $\mathscr{X}$ is a quasi-simple \D-scheme then all the restriction maps of the structure sheaf induce isomorphisms on the constant rings (see Proposition \Ref{Cst_faisc_cst_quasi_simple}). The next section is devoted to the computation of the coarse space of leaves for quasi-simple schemes. We prove (see Theorem \Ref{Existence_Coarse_Sp_Leaves_Qu_Simple}):

\smallskip

\begin{theoremeintro}
Let $\mathscr{X}$ be a quasi-simple \D-scheme. Then, $\mathscr{X}$ has a coarse space of leaves  $t_{\mathscr{X}}: \mathscr{X}\longto C(\mathscr{X})$. Furthermore, $C(\mathscr{X})$ is the spectrum of a field.
\end{theoremeintro}

\noindent In order to prove this theorem, we prove a lemma (Lemma \Ref{lemme_hyper_technique}) on the patching of universal objects. Section 5 is devoted to the 
affine statements of our geometrization of Lemma S-vdP. In Section 6, we prove the main result of this paper. The proof consists in showing that one can reduce successively to the case where the total space $X$ is affine and to the case where both $X$ and $S$ are affine. Finally, in an appendix, we give an overview of commutative \D-algebra. The principal difference with the classical references \Cite{kaplansky} or \Cite{DiffAlgGroupsKolch} is that we deal with rings endowed with any family of derivations (instead of rings with one derivation or with a finite family of commuting derivations).  The main content of this appendix is the two proofs of the primalily of $\mathfrak{p}_\#$ when $\mathfrak{p}$ is a prime ideal of some $\Q$-\D-algebra. This section contains also a subpart dedicated to simple \D-rings. Among other results, we prove that simple \D-rings are always irreducible, without any condition on the characteristic or the reducedness. It relies of the following lemma (see Lemma \Ref{Super_lemme}).

\begin{lemmaintro}
Let $R$ be a ring. Let $I$ and $J$ be two ideals of $R$ such that $I \cap J=0$. Let $x\in I$ and $y\in J$. Then, 
\[\fa \theta, \theta'\in \widehat{\Theta}(R), \quad\theta(x)\cdot \theta'(y)=0.\]
\end{lemmaintro}

\noindent A geometric corollary of it (see Corollary \Ref{Coro_geo_irr}) says that if a scheme $X$ can be decomposed as $X=F_1\cup F_2$ where the $F_i$'s are closed sets, then the same decomposition holds when replacing  each $F_i$ by 
\[ [F_i]:= \bigcap_{\chv{V} \text{ vec. field}}\text{invariant subset of $F_i$ under $\chv{V}$}.\]
We have also included a paragraph on the colon ideals in differential rings.

\bigskip

The author wishes to thank Fran\c{c}ois Ollivier  for helpful and interesting discussions.

\section{Vector fields, leaves and trajectory for schemes}

In this section, we recall the definitions and main properties of vector fields, leaves and of the trajectory, as introduced and studied in our previous paper \Cite{Traj1}. We refer to this paper for examples and further comments.  We also show that one can generalize these constructions for schemes endowed with several vector fields. It is to be noted that we don't make the usual assumption that the derivations commute. These generalizations are mostly straightforward, excepted for Theorem \Ref{theo_traj}, which gives the definition of $\Traj{V}(x)$. The proof of this theorem is postponed until the appendix, see Propositions  \Ref{P_diese_premier_bis} and \Ref{primality_ideals}. Let us start this paper with some notations and conventions on \D-rings.

\bigskip

\subsection{\D-rings}\label{D_Rings}

In all this paper, $\Delta$ will be a fixed set. We  denote by $d$ the cardinal of $\Delta$ if finite and $\infty$ otherwise. The elements of $\Delta$ will be denoted by $\partial, \partial', \partial_i$, etc. When $\Delta$ is finite, we will denote $\Delta:= \{\partial_1, \ldots, \partial_d\}$. We will work with rings endowed with a familiy of  derivations indexed by $\Delta$ or, equivalently, with schemes endowed with a family of vector fields (see Definition \Ref{Def_Vect_Fields}). 
 More precisely\footnote{All rings are commutative and with unit.}, we follow the definition of \Cite{buium}:
\begin{definition}
A \emph{\D-ring} is a ring $R$ with a map 
\[ 
 \phi : \fonctionb{\Delta}{\mathrm{Der}(R)}{\partial}{\phi(\partial)}.
 \]
When $d\neq \infty$ and when for all $\partial, \partial'\in \Delta$ one has $[\phi(\partial), \phi(\partial')]=0$, we say that $(R, \phi)$ is a \emph{partial \D-ring}. When $d=1$, it is called an \emph{ordinary \D-ring}.
\smallskip

A morphism $\varphi: (R, \phi)\longto(R', \phi')$ between two \D-rings is a ring morphism $\varphi:R\longto R'$ such that  the diagram
\[ 
\xymatrix{
R \ar[d]_{\phi(\partial)} \ar[r]^\varphi & R' \ar[d]^{\phi'(\partial)}\\
R \ar[r]^\varphi & R'
}
 \]
 commutes for all $\partial\in \Delta$.  

\smallskip 
 
 The category of \D-rings will be denoted by $\DRng$. 
 
\smallskip 
 
 Similarly, if $k$ is a \D-ring, one has the category $\DAlg{k}$ of $k$-\D-algebras.
\end{definition} 
 \bigskip
 
 If $(R, \phi)$ is a \D-ring and if there is no danger of confusion, we will write $\partial$ insted of $\phi(\partial)$.
   Given a \D-ring $R$, we will denote its ring of constants by $C(R)$ (though we could also have chosen the notation $R^{\Delta}$). 
 Now, a few words on \emph{differential operators}:
 \begin{itemize}
\item $\Theta(R)$: given a \D-ring $R$,  it will denote the free monoid generated by $\Delta$. (Thus, it does only depend on $\Delta$). Every $\theta\in \Theta(R)$ can be uniquely written
\[ 
\theta=\partial_{1}\partial_{2}\cdots \partial_{n},
 \]
 with the $\partial_i$'s in  $\Delta$. The \emph{order of $\theta$}, denoted by $e(\theta)$, is  the integer $n$. 
 \item $\Theta^{\textit{ab}}(R)$: given a \emph{partial} \D-ring $R$, it will denote the free \emph{commutative} monoid generated by  $\Delta$. Every $\theta\in \Theta^{\textit{ab}}(R)$ can be uniquely written
\[ \theta=\prod_{i=1}^{d} {\partial_i}^{e_i(\theta)}
 \]
In this case,  the order of $\theta$ is the integer $e(\theta)=\sum_i e
 _i(\theta)$. 
 \item $\widehat{\Theta}(R)$: given a  ring $R$, it will denote the free monoid generated by  $\mathrm{Der}(R)$. Every $\theta\in \widehat{\Theta}(R)$ can be uniquely written
\[ 
\theta=\delta_1\delta_2\cdots \delta_n,
 \]
 with the $\delta_i$'s in  $\mathrm{Der}(R)$. The order of $\theta$ is $e(\theta):=n$. 
\end{itemize}
Given $R$ a ring (resp. a \D-ring, a partial \D-ring), $\widehat{\Theta}(R)$ (resp. $\Theta(R)$, $\Theta^{\textit{ab}}(R)$) acts in a natural way on $R$. A \emph{\D-ideal} $I$ of a \D-ring $R$ is an ideal such for all $\theta\in \Theta(R)$, $\theta(I)\subset I$. Equivalently,  $I$ is \D-ideal iff $\partial(I)\subset I$ for all $\partial\in \Delta$. For any set $I\subset R$, there exists a unique smallest \D-ideal containing $I$, that will be denoted by $\langle I \rangle$.
\smallskip

The appendix \Ref{appendix} is devoted to commutative \D-algebra. For further reference on differential algebra in the partial and in the ordinary case, see \Cite{kaplansky}, \Cite{singvdp} or \Cite{DiffAlgGroupsKolch}.

\subsection{Vector fields and \D-schemes} Vector fields are the exact analogues for schemes to derivations for rings. We define:

\begin{definition}\label{Def_Vect_Fields}
Let $X$ be a scheme. A \emph{vector field on X} is a derivation of the structure sheaf $\mathcal{O}_X$. 
\end{definition}

 In the case of a smooth manifold $M$, it is well known that there is a one-to-one correspondence between vector fields defined as global sections of the tangent bundle $TM \longto M$ and derivations of the smooth struture sheaf $\mathcal{O}_M$.        
 For schemes, one also has such a correspondence. In \Cite{EGA44} \S \textbf{(16.5)}, given a $S$-scheme $X$, Grothendieck defines the tangent bundle of $X/S$. It is a $S$-scheme, denoted by $T_{X/S}$ , with a $S$-morphism to $X$:
\[ \xymatrix@R=0.6cm{T_{X/S} \ar[d]^{\pi}\\X} \]
He proves that the $S$-section of $\pi$ correspond to the $\mathcal{O}_S$-derivations of $\mathcal{O}_X$. So,
in the case where $X$ is viewed as a $\Z$-scheme, one gets a correspondence between
the sections of $\pi : TX \longto X$ and the group of vector fields of $X$. The $\mathcal{O}_X$-module
of $S$-sections of $\pi$ is the dual of $\Omega^1_{X/S}$. We will denote it by $\mathscr{T}_{X/S}$ (or by $\mathscr{T}_X$ when $ S = \Sp \Z$). Hence, $\mathscr{T}_X(X)$ is the group of vector fields on $X$.

\begin{fact}
$\mathscr{T}_X(X)$ is an $\mathcal{O}_X(X)$-module and a Lie algebra.
\end{fact}

\pproof{Given $\boldsymbol{\partial}^a:=(\partial^{\,a}_U)_U$ and  $\boldsymbol{\partial}^b:=(\partial^{\,b}_U)_U$ two vector fields on $X$, ones defines $[\boldsymbol{\partial}^a, \boldsymbol{\partial}^b]$ to be the derivation of $\mathcal{O}_X$ whose action on $\mathcal{O}_X(U)$ is given by the derivation $\left [\partial^{\,a}_U, \partial^{\,b}_U\right ]$. The structure of $\mathcal{O}_X(X)$-module is clear. }
\smallskip

Given a scheme $X$ and vector field $\chv{V}$ on $X$, one can associate to each  point $x\in X$ a Zariski tangent vector, that we will denote by $\chv{V}(x)$ (see \Cite{Traj1} for the construction). The application
\[ \fonctionb
{\mathscr{T}_X(X) }
{T_x X}
{\chv{V}}
{\chv{V}(x)} \]
is $\mathcal{O}_X(X)$-linear. 
\smallskip

We say that two vector fields $\chv{V}$ and $\chv{W}$ \emph{commute} if $[ \chv{V}, \chv{W}]=0$.

\begin{definition}\label{DeltaScheme}
A \emph{\D-scheme} is a scheme $X$ endowed with a family $\fchv{V}=(\chv{V}_\partial)_{\partial\in \Delta}$ of vector fields on $X$. \smallskip

 Given two \D-schemes $\mathscr{X}=(X, \fchv{V})$ and $\mathscr{Y}=(Y, \fchv{W})$ a morphism $f : \mathscr{X} \longto \mathscr{Y}$ is a morphism of schemes $f=(\varphi, \theta):X \longto Y$ such that for all $\partial\in \Delta$ the diagram
\[ 
\xymatrix@R=8mm
{
\ar[d]^{\chv{W}_\partial} \mathcal{O}_Y \ar[r]^-{\theta} &\varphi_* \mathcal{O}_X \ar[d]^{\varphi_* \chv{V}_\partial} \\
\mathcal{O}_Y \ar[r]^-{\theta}& \varphi_* \mathcal{O}_X 
}
 \]
 commutes. 
\smallskip 
 
 The category of \D-schemes is denoted by $\DSch$.
\end{definition}

\smallskip

\D-schemes will be denoted by letters $\mathscr{X}$, $\mathscr{Y}$, etc. and the underlying schemes by $X$, $Y$, etc. The corresponding families of vector fields will be denoted by $\fchv{V}$, $\fchv{W}$, etc. Instead of writting ``for $\chv{V}_\partial$, with $\partial\in \Delta$'', we will simply write, if there is no danger of misunderstanding, ``for $\chv{V}\in \Delta$''.  If $(P)$ is a property for schemes (affine, reduced, irreducible, of characteristic 0, etc.) and if $\mathscr{X}$ is a \D-scheme, we will say that $\mathscr{X}$ has $(P)$ if $X$ has.  Given a \D-scheme $\mathscr{X}$, the structure sheaf of $X$ is naturally a sheaf of \D-rings. We will still denote it by $\mathcal{O}_X$.

\smallskip
If $f$ is a morphism of \D-schemes and $x\in X$, then $f$ sends the tangent vectors $\fchv{V}(x)$ to  $\fchv{W}(f(x))$:

\begin{proposition}
Let $\mathscr{X}=(X, \fchv{V})$ and $\mathscr{Y}=(Y,\fchv{W})$ be two \D-schemes and $f: \mathscr{X}\longto \mathscr{Y}$ be a morphism. Let $x\in X$. Then, for all $\partial\in \Delta$,
$$T_{x} f \bullet \chv{V}_\partial (x)=i_{x}\circ \chv{W}_\partial (f(x))$$
where $i_{x} : \kappa(f(x))\longto \kappa(x)$ is the inclusion of residual fields induced by $f$.
\end{proposition}

\pproof{See Proposition 2.3 of \Cite{Traj1}. }

\bigskip

There is a functor $(\DRng)^{\mathrm{op}} \longto \DSch$ that maps a \D-ring $(R, \phi)$ to the scheme $\Sp R$ endowed with the vector fields associated to the derivations $\phi(\partial)$. We denote it by $\DSp$.  This functor $\DSp$ is a right adjoint to the functor $\mathcal{O}(-) : \DSch \longto (\DRng)^{\mathrm{op}}$ of global sections. 

\smallskip

All fibered products exist in   $\DSch$. Given \D-schemes $\mathscr{X}, \mathscr{Y}$ above $\mathscr{S}$, the underlying scheme of $\mathscr{X}\times_{\mathscr{S}}\mathscr{Y}$ is $X\times_S Y$. Equivalently, cartesian squares (and more generally, all limits) commute with the forgetfull functor $\DSch\longto \Sch$. This is because this functor has a left adjoint (the affine version of this statement is proved in \Cite{gillet}, for instance).

\subsection{Leaves and trajectory} First, let us define the leaves.

\begin{definition}\label{Defi_leaf}
Let $\mathscr{X}=(X, \fchv{V})$ be a \D-scheme. Let $\eta \in X$. We say that \emph{$\eta$ is a leaf of $\mathscr{X}$} (or \emph{a leaf for $\fchv{V}$}) if $\chv{V}(\eta)=0$ for all $\chv{V}\in \Delta$. The set of leaves of $\mathscr{X}$ will be denoted by $X^{\fchv{V}}$.
\end{definition}

\noindent If $R$ is a \D-ring, then the leaves of $\DSp R$ are exactly the prime \D-ideals of $R$ (see \Cite{Traj1}). So, the sets $X^{\fchv{V}}$ are generalizations of the more common differential spectrum ``$\mathrm{diff\,Spec}$'' (see \Cite{KovacicDiffSchemes}). Intuitively, the leaves of $\mathscr{X}$ correspond to (the generic points of) the irreducible closed subvarieties of $X$ tangent to the vector fields $\chv{V}_i$. 
For instance, one has

\begin{proposition}
Let $X$ be a scheme of characteristic zero and let $C$ be an irreducible component of $X$. Then, the generic point $\eta_C$ of $C$ is a leaf for every vector field. 
\end{proposition}

\pproof{
It follows from \emph{(ii)} of Proposition \Ref{primality_ideals}. 
}

\bigskip

\noindent The following theorem defines the trajectory and gives some properties of it:

\bigskip

\begin{theorem}\label{theo_traj}
Let  $\mathscr{X}=(X, \fchv{V})$ be a $\Q$-\D-scheme and let $x\in X$.
\begin{itemize}

\item[(i)]The set $\left \{ \eta \in X\,\,\middle | \,\,\eta \leadsto x \aand \eta \in X^{\fchv{V}}\right \}$ has a least element. We call it \emph{the trajectory of $x$ (under the action of $\fchv{V}$)} and denote it by $\Traj{V}(x)$.

\item[(ii)]The map $\Traj{V} : X \longto X^{\fchv{V}}$ is continuous and open for the induced topo\-logy on $X^{\fchv{V}}$.

\item[(iii)]If $(Y, \fchv{W})$ is another  $\Q$-\D-scheme and $f : \mathscr{X} \longto \mathscr{Y}$ is a morphism, then
$
f(\Traj{V}(x))=\Traj{W}(f(x)).
$
\end{itemize}
\end{theorem}

\bigskip

\pproof{\emph{(i)}: as in the case $d=1$ (see Theorem 2.5 of \Cite{Traj1}), it suffices to check the result for affine schemes.  So, one has to check that if $R$ is a $\Q$-\D-algebra  and if 
$\mathfrak{p}\in \Sp R$, then the \D-ideal
 \[\mathfrak{p}_{\#} :=\left\lbrace f\in R \tq \fa  \theta \in \Theta(R), \quad \theta(f)\in \mathfrak{p} \right\rbrace   \]
 is prime. We give two proofs of this result, in Proposition \Ref{primality_ideals} and \Ref{P_diese_premier_bis}.  For \emph{(ii)} and  \emph{(iii)}, see Proposition 2.6 and 2.7 of \Cite{Traj1}. }

\bigskip

Let us mention here an easy result on the specialization order. We will use it implicitely in what follows, and especially in Section 6.

\begin{fact}
Let $X$ and $Y$ be two topological spaces, let $f:X \longto Y$ be a continuous function. Then,
\[
\fa \eta, x\in X, \qquad \eta\leadsto x\ \impl\ f(\eta)\leadsto f(x).
\]
\end{fact}

\smallskip

\pproof{
Let $x,\eta\in X$ and assume $\eta\leadsto x$. Let us prove that $f(x)\in\overline{ \{f(\eta)\}}$. Let $F$ be closed set of $Y$ such that $f(\eta)\in F$. Then, $\eta\in f^{-1}(F)$. Since $x\in \overline{\{\eta\}}$, one has $x\in  f^{-1}(F)$ and so $f(x)\in F$.
}

\bigskip

\section{Simple and quasi-simple \D-schemes}

In this section, we introduce simple and quasi-simple \D-schemes. Simple \D-schemes will replace the base field $k$ of Lemma S-vdP in the geometrized statement. Actually, in a way, simple \D-schemes stand for, in the differential setting, the points $\Sp K$ (with $K$ a field) of the classical setting. The quasi-simple \D-schemes are a slight generalization of simple \D-schemes. We will prove for them the existence of a coarse space of leaves. One can show that the coarse space of leaves, when it exists, is not always a good approximation of what one would hope --- quasi-simple \D-schemes can provide such examples.

\subsection{Simple \D-schemes}

A \D-ring $R$ is said to be \emph{simple} (see Appendix \RefPar{section_simple_D_rings}) if the only \D-ideals of $R$ are $0$ and $R$. Fact \Ref{char_simple} tells us that a $\Q$-\D-ring is simple iff
the only  prime \D-ideal of $R$ is $0$. This leads us to the following definition:

\begin{definition}\label{DScheme_simple}
A \D-scheme $\mathscr{X}$ is said to be \emph{simple} if it  is irreducible, and   if its only leaf is its generic point.
\end{definition}

If $k$ is any field of characteristic $0$ and if $k[x]$ is endowed with the derivation $x'=1$, then $k[x]$ is simple. So, $\A{1}{k}$ with the vector field $\partial/\partial x$ is simple.

\begin{proposition} \

\begin{itemize}\setlength{\itemsep}{0mm}
\item[(i)] Let $R$ be a \D-ring. Then, $
R\text{ simple }\impl \DSp R \text{ simple}$. The converse is true when $R$ does not have  zero divisors.
\item[(ii)] A \D-scheme $\mathscr{X}$ is simple iff it is irreducible and all its open affine \D-subschemes are simple.
\item[(iii)] A \D-scheme $\mathscr{X}$ is simple iff it is irreducible and it has an affine atlas $(\DSp R_i)_i$ where the  its open affine \D-subschemes are simple.
\item[(iv)] Let $\mathscr{X}$ be a simple $\Q$-\D-algebra. Then, $\mathscr{X}$ is integral.
\end{itemize}
\end{proposition}

\pproof{
For \emph{(i)} and \emph{(iv)}, use Fact \Ref{char_simple} and Proposition \Ref{Properties_simple}.\emph{(iii)}. For \emph{(ii)} and \emph{(iii)}, use the fact that a point is a leaf iff it is a leaf in some open affine neighborhood. }

\subsection{The constant sheaf of a simple \D-scheme}

We would like to prove an analogue for simple \D-schemes of Proposition \Ref{Properties_simple}.\emph{(ii)}, which says that $C(R)$ is a field when $R$ is a simple \D-field. We want to associate to any simple \D-scheme a scheme which would be the analogue of the ring of constants. In a way, such a scheme would parametrize the leaves of $\mathscr{X}$ --- we will come back to this idea in section \ref{Sect_coarse_space_leaves}. In this subsection, we prove that given a simple \D-scheme $\mathscr{X}$, the sheaf of constants of $\mathscr{X}$ is constant, with a field as stalks. It will be an important ingredient in the construction of a coarse space of leaves $C(\mathscr{X})$ for (quasi-)simple \D-schemes. 
\smallskip

Given a topological space $X$ endowed with a sheaf of \D-rings $\mathscr{F}$, we will denote by $\mathscr{F}^{\Delta}$ and call the \emph{the constant sheaf associated to $\mathscr{F}$} the sheaf defined by
$
\mathscr{F}^{\Delta}(U):=C (\mathscr{F}(U))
$
for all open set $U$. In the case where $\mathscr{F}$ is the structure sheaf of a \D-scheme $\mathscr{X}$, we will call $\mathcal{O}_X^{\Delta}$ the \emph{constant sheaf of $\mathscr{X}$}. So, the main result of this paragraph is:

\begin{proposition}\label{Cst_sh_cst}
Let $\mathscr{X}$ be a simple \D-scheme. 

\begin{itemize}
\item[(i)] For every nonempty open sets $ U\subset V$, the constant trace of the restriction map
$ 
\rho_{V\to U}^\Delta : \mathcal{O}_X^\Delta(V) \longto \mathcal{O}_X^\Delta(U)
 $
 is an isomorphism.
 \item[(ii)] Furthermore, for all nonempty open set $U$, the ring  $\mathcal{O}_X^\Delta(U)$  is a field.
\end{itemize}
\end{proposition}

We will use the following lemma, whose proof is left to the reader.

\begin{lemma}\label{iso_base_sh}
Let $X$ be an irreducible topological space and let $(V_i)_{i\in I}$ be a basis of nonempty open sets for $X$. Let $\mathscr{F}$ be a sheaf of abelian groups on $X$ verifying
\[ 
\fa (i,j)\in I^{2}, \quad V_i \subset V_j \ \impl\  \rho_{V_j\to V_i} : \mathscr{F}(V_j)\to \mathscr{F}(V_i) \text{ is an isomorphism}.
 \]
 Then, for all nonempty open sets $U$ and $V$ such that $U \subset V$, the map \[ \rho_{V\to U}: \mathscr{F}(V)\longto \mathscr{F}(U) \text{ is an isomorphism.} \] 
\end{lemma}

\smallskip

\pproofbis{of Proposition \Ref{Cst_sh_cst}}{
First, let us prove the result in the affine case. Let $R$ be a simple \D-algebra and $\mathscr{X}:= \DSp R$. A basis of nonempty open sets is given by the $D(f)$, for $f\in R\setminus \Nil(R)$. Let $f$ and $g$ such that $D(f)\subset D(g)$. Then, in the commutative diagram
\[ 
\boite{ \xymatrix@C=4mm@R=7mm{
&  \mathcal{O}_{X}^{\Delta}\left ( X \right ) \ar[rd]^{\rho_{2}}\ar[ld]_{\rho_{1}}\\
  \mathcal{O}_{X}^{\Delta}\left ( D(g) \right )  \ar[rr]_{\rho_{3}}&& \mathcal{O}_{X}^{\Delta}\left ( D(f) \right )
 }},
 \]
 whose in fact can be written more explicitely
 \[ 
\boite{ \xymatrix@C=4mm@R=7mm{
& C(R) \ar[ld]_{\rho_{1}}\ar[rd]^{\rho_{2}}\\
  C(R_g) \ar[rr]_{\rho_{3}}&& C(R_f)
 }},
 \]
 the maps $\rho_1$ and $\rho_2$ are isomorphisms, by Proposition \Ref{Properties_simple}.\emph{(iv)}. Hence, $\rho_3$ is also an isomorphism, and the result follows from Lemma \Ref{iso_base_sh}.
 
 \smallskip
 
 Now, if $\mathscr{X}$ is any \D-scheme, let us consider the basis of nonempty affine open sets. Let $U$ and $V$ be two such sets, with $U\subset V$. By the previous case, one knows that in the commutative diagram
 \[ 
\boite{
\xymatrix@C=4mm@R=7mm
{
\mathcal{O}_X^{\Delta}(V) \ar[rr]^-{\rho_1}\ar[rd]_-{\rho_2}
&&
\mathcal{O}_X^{\Delta}(U) \ar[ld]^-{\rho_3} 
\\
&
\mathcal{O}_X^{\Delta}(U \cap V)  
}
} 
  \]
the maps $\rho_2$ and $\rho_3$ are isomorphisms, and so is $\rho_1$. The same lemma gives us the conclusion.
}

\subsection{Quasi-simple \D-rings and quasi-simple \D-schemes}

In this subsection, we define quasi-simple \D-rings and \D-schemes. Quasi-simple \D-schemes have the same formal properties as simple \D-schemes, regarding the existence of a coarse space of leaves (see Theorem \Ref{Existence_Coarse_Sp_Leaves_Qu_Simple}). Let us start by an example.

\smallskip

In this example, $d=1$. Let $k$ be a field, endowed with the zero derivation. We consider the $k$-\D-algebra $R$ defined by 
\[ R:= \left \{ \begin{array}{l}
k[x]\\
x'=x
\end{array}\right. . \]
We call $\DSp R$ the \emph{affine line endowed with the radial field}. This \D-scheme has two leaves: the generic point $\eta$ and the point $0$. The \D-ring $R$ satifies:

\begin{fact}
Assume $k$ to be of characteristic zero. Then, $C(\Frac R)=C(R)$.
\end{fact}

\pproof{First, it is easy to check that $C(R)=k$. Then, let  $f,g\in k[x]$, with $g\neq  0$. Denote $h:=f/g$, and assume 
$
h'=\left ( {f}/{g} \right )'=0.
$
This implies $f'g-g'f=0$. In other words,
$$
\frac{f}{g}=\frac{f'}{g'}.
$$
Hence, one has:
$$
\frac{f}{g}=\frac{f'-\deg(f)\cdot f}{g'-\deg(f)\cdot g}.
$$
This identity is valid only if  $g'-\deg(f) g \neq 0$. But, in the case when $g'=\deg(f) g$, the identity $f'g-g'f=0$ can be written
$(f'-\deg(f)f)g=0$. Denoting $n=\deg(f)$, one then has
$
g= \lambda x^n$ and $ f=\mu x^n
$
and so $h\in k$.
\smallskip

When the transformation of the denominator if possible, then, by iteration, one can reduce the numerator of $h$ to a constant: there exist $\lambda\in $ and $p\in k[x]$ such that $h=\lambda/p$. It is easy to check then that $p\in k$.
}

\smallskip

This leads us to the following:

\begin{definition}
A \D-ring $R$ is said to be \emph{quasi-simple} if $R$ is a domain and if the map $C(i) : C(R)\longto C(\Frac R)$ is an isomorphism.\smallskip

A \D-scheme $\mathscr{X}$ is said to be \emph{quasi-simple} if it is irreducible and if it has an atlas $(\DSp R_i)_{i\in I}$ where the \D-rings $R_i$ are all quasi-simple.
\end{definition}

Of course, by Proposition \Ref{Properties_simple}.\emph{(iv)}, simple \D-rings (resp. schemes) are quasi-simple \D-rings (resp. schemes). Now, our aim is to prove that quasi-simple \D-schemes verify the same property as simple ones: their sheaf of constants is constant.

\smallskip

\begin{proposition}\label{Cst_faisc_cst_quasi_simple}
Let $\mathscr{X}$ be a quasi-simple \D-scheme. 

\begin{itemize}\setlength{\itemsep}{0mm}
\item[(i)] For every nonempty open sets $ U\subset V$, the constant trace of the restriction map
$
\rho_{V\to U}^\Delta : \mathcal{O}_X^\Delta(V) \longto \mathcal{O}_X^\Delta(U)
$
 is an isomorphism.
 \item[(ii)] Furthermore, for all nonempty open set $U$, the ring  $\mathcal{O}_X^\Delta(U)$  is a field.
\end{itemize}
\end{proposition}
\smallskip

We will use the two following lemmas. The first lemma is a purely geometric result on integral schemes.

 \begin{lemma}\label{lemma_classique_integre}
Let $X$ be an integral scheme. Let $U$ and $V$ be two nonempty open sets of $X$ with $U\subset V$. Then,

\begin{itemize}\setlength{\itemsep}{0mm}
\item[a)] The rings $\mathcal{O}_X(U)$ and $\mathcal{O}_X(V)$ do not have zero divisors.
\item[b)] The morphism $\mathcal{O}_X(V)\longto\mathcal{O}_X(U)$ is injective.
\item[c)] If moreover $V$ is affine then  $\Frac \mathcal{O}_X(V)\longto\Frac\mathcal{O}_X(U)$ is an isomorphism.
\end{itemize}
\end{lemma}

\pproof{
The point \emph{a)} is classical, see \Cite[II.3]{hartshorne}. Let us prove \emph{b)}. We first discuss the case when $V$ is affine. Hence, let $A$ be domain and $U$ a nonempty open set of $X:=\Sp A$. Let $f\in A$ such that $\varnothing\neq D(f)\subset U$. The diagram
\[ 
\xymatrix@R=6mm@C=3mm{
& \mathcal{O}_{X}\left ( X \right )\ar[ld]_-{\rho_{1}} \ar[rd]^-{\rho_{2}}\\
 \mathcal{O}_{X}\left ( U \right )\ar[rr]_-{\rho_{3}}&& \mathcal{O}_{X}\left ( D\left ( f \right ) \right )
}
 \]
 is commutative and can be written more explicitely
 \[ 
\boite{\xymatrix@R=6mm@C=3mm{
&A\ar[ld]_-{\rho_{1}} \ar[rd]^-{\rho_{2}}\\
\mathcal{O}_{X}\left ( U \right )\ar[rr]_-{\rho_{3}}&&A_{f}}}. 
  \]
  But, $\rho_2$ is injective and thus so is $\rho_1$. Let us discuss now the general case. Let $X$ be integral scheme and $U\subset V$ two nonempty open sets. We want to show that 
  \[\varphi : \mathcal{O}_X(V)\longto \mathcal{O}_X (U)   \]
  is injective. Let $x\in \mathcal{O}_X(V)$ such that $\varphi(x)=0$. For all affine open set $W$ of $V$, the restriction $\rho_W : \mathcal{O}_X(W)\longto \mathcal{O}_X (U\cap W)$ is injective, by the previous case. So, $x_{\mid W}=0$, and so is $x$.

 \smallskip
 
 Let us prove \emph{c)} now. It is sufficient to prove that if $R$ is a domain and if $U$ is a nonempty open set, then $\Frac A \longto \Frac\mathcal{O}_X(U)$, which is well-defined by \emph{b)}, is an isomorphism. Let $f\in R$ such that $\varnothing \neq D(f) \subset U$. The diagram
 \[ 
\xymatrix@R=4mm@C=4mm
{
\ar[rd]_-{\rho_2}\Frac R \ar[rr]^-{\rho_1} &&\Frac \mathcal{O}_X(U)\ar[ld]^-{\rho_3}\\
&\Frac A_f
} 
  \]
  commutes. Since $\rho_2$ is bijective, $\rho_3$ is surjective. But, $\rho_3$ is also injective and hence, $\rho_3$ is bijective. So, $\rho_1$ is bijective.
}

\smallskip

\begin{lemma}\label{lemma_affine_subset_affine_quasi_simple}
Let $R$ be a  quasi-simple \D-ring and $X:=\DSp R$. Let $U$ be a nonempty open set of $X$. Then, the \D-ring $\mathcal{O}_X(U)$ is quasi-simple.
\end{lemma}

\smallskip

\pproof{
The diagram
\[ 
\xymatrix@C=2cm{
R \ar[r]^-{\varphi}\ar[d]^-{i_{1}}& \mathcal{O}_{X}\left ( U \right )\ar[d]_-{i_{2}}\\
\Frac R \ar[r]^-{\Frac \varphi} & \Frac \mathcal{O}_{X}\left ( U \right )
}
 \]
 commutes and so does its trace on constants
 \[ 
\boite{\xymatrix@C=2cm{
C(R) \ar[r]^-{C(\varphi)}\ar[d]^-{C({i_{1}})}&C( \mathcal{O}_{X}\left ( U \right ))\ar[d]_-{C({i_{2}})}\\
C(\Frac R) \ar[r]^-{C({\Frac \varphi})} & C({\Frac \mathcal{O}_{X}\left ( U \right )}
)}}. 
  \]
By assumption, $C(i_1)$ is an isomorphism. By \emph{c)} of Lemma \Ref{lemma_classique_integre}, $\Frac \varphi$ is an isomorphism and so is $C(\Frac \varphi)$. Hence, $C(i_2)$ is surjective, but is also injective, as $i_2$ is. So, $C(i_2)$ is an isomorphism, and $ \mathcal{O}_{X}\left ( U \right )$ is quasi-simple.
}

\bigskip

\pproofbis{of Proposition \Ref{Cst_faisc_cst_quasi_simple}}
{
Let us fix the notations. $\mathscr{X}$ is a quasi-simple \D-scheme, $(U_i)_{i\in I}$ is an open covering of $X$ by affine schemes, such that $\mathcal{O}_X(U_i)$ is quasi-simple, for all $i$. Let $(V_\alpha)_{\alpha\in A}$ be the basis of all nonempty open affine sets included in a least one $U_i$. Let $V_\alpha$ and $V_\beta$ two elements of this basis, such that $V_\alpha\subset V_\beta$. By lemma \Ref{lemma_affine_subset_affine_quasi_simple}, one knows that the \D-rings $\mathcal{O}_X(V_\alpha)$ and $\mathcal{O}_X(V_\beta)$ are quasi-simple. Hence, in the diagram
\[ 
\xymatrix@C=2cm{
C({\mathcal{O}_{X}\left ( V_{\beta} \right )}) \ar[r]^-{C({\varphi})}\ar[d]^-{C({i_{1}})}&C({ \mathcal{O}_{X}\left ( V_{\alpha} \right )})\ar[d]_-{C({i_{2}})}\\
C({\Frac \mathcal{O}_{X}\left ( V_{\beta} \right )}) \ar[r]^-{C({\Frac \varphi})} & C({\Frac \mathcal{O}_{X}\left ( V_{\alpha} \right )})
},
 \]
 $C(i_1)$ and $C(i_2)$ are isomorphisms. Then, since $V_\beta$ is affine and $X$ integral, by Lemma \Ref{lemma_classique_integre}, $C(\Frac \varphi)$ is also an isomorphism. So, $C({\varphi})$ is an isomorphism. Now, one can conclude by Lemma \Ref{iso_base_sh}.
}

\bigskip

As a first application, let us prove the following:

\begin{corollary}\label{Carac_Sch_Quasi_Simple} \
\begin{itemize}\setlength{\itemsep}{0mm}
\item[(i)] $\mathscr{X}$ is quasi-simple iff for all nonempty open affine set $U$ of $X$, the \D-ring $\mathcal{O}_X(U)$ is quasi-simple.
\item[(ii)] Any nonempty open subscheme of a quasi-simple \D-scheme is  quasi-simple.
\end{itemize}
\end{corollary}

\pproof{The point \emph{(ii}) is an easy consequence of \emph{(i)}. Let us prove the latter.
The direction $\Longleftarrow$ is clear. Conversely, let us denote by $(V_i)_{i}$ a covering of $X$ by open affines sets such that the \D-rings $\mathcal{O}_X(V_i)$ are quasi-simple.  Let $U$ be a nonempty affine open set of $X$. Let $U_{i_0}$ be one of the element of the covering. Since $X$ is irreducible,  $U\cap U_{i_{0}}\neq \varnothing$. As previously, the diagram
\[ 
\xymatrix@C=2cm{
C({\mathcal{O}_{X}\left ( U  \right )}) \ar[r]^-{C({\varphi})}\ar[d]^-{C({i_{1}})}&C({ \mathcal{O}_{X}\left ( U\cap U_{i_{0}} \right )})\ar[d]_-{C({i_{2}})}\\
C({\Frac \mathcal{O}_{X}\left ( U \right )}) \ar[r]^-{C({\Frac \varphi})} & C({\Frac \mathcal{O}_{X}\left ( U\cap U_{i_{0}} \right )})
}
 \]
 commutes. Now, by Lemma \Ref{lemma_affine_subset_affine_quasi_simple}, $\mathcal{O}_X(U\cap U_{i_{0}})$ is quasi-simple, and so $C(i_2)$ is an isomorphism. Since $U$ is affine, by Lemma \Ref{lemma_classique_integre}.\emph{c)}, $C(\Frac \varphi)$ is an isomorphism. Last, by Proposition \Ref{Cst_faisc_cst_quasi_simple}, the constant trace $C(\varphi)$ of the restriction is an isomorphism. Hence, $C(i_1)$ is an isomorphisme, \ie $\mathcal{O}_X(U)$ is quasi-simple.
}

\section{The coarse space of leaves: case of quasi-simple \D-schemes}\label{Sect_coarse_space_leaves}

\subsection{The coarse space of leaves}

Let $\mathscr{X}$ be a \D-scheme of characteristic zero. We would like to define a ``space '' $T$ that classify the leaves of $\mathscr{X}$, a space whose points would intuitively correspond to the trajectories of $\mathscr{X}$. Such a space would be endowed with a map $\varphi : X \longto T$ that should verify
\begin{equation}\label{comp_traj}
\fa (x,y)\in X^2 \qquad \Traj{V}(x)=\Traj{V}(y)\quad \impl \quad \varphi(x)=\varphi(y).
\end{equation}
Now, let $T$ be a $\Q$-scheme, endowed with the zero vector fields. Let us consider $\varphi : \mathscr{X}\longto (T, \vec{\mathbf{0}})$ a morphism of \D-schemes and let us verify that (\ref{comp_traj}) stands for $\varphi$. Let $x,y\in X$. By Theorem \Ref{theo_traj}, one has
\[ 
\mathrm{Traj}_{\vec{\mathbf{0}}}(\varphi(x))=\varphi(x)=\varphi(\Traj{V}(x)).
 \]
 Hence, indeed, if $ \Traj{V}(x)=\Traj{V}(y)$ then $\varphi(x)=\varphi(y)$. This leads us to the definition:
 
 \smallskip
 
 \begin{definition}\label{Def_coarse_space_of_leaves}
Let $\mathscr{X}$ be a \D-scheme of characteristic zero. A \emph{coarse space of leaves for $\mathscr{X}$} will be a $\Q$-scheme $T$ with a map $\varphi : \mathscr{X}\longto (T, \vec{\mathbf{0}})$ universal for this property.
\end{definition}

\smallskip

This means that  every morphism $\mathscr{X}\longto (T', \vec{\mathbf{0}})$ factors uniquely through $\varphi$. 
\emph{A priori}, such a coarse space of leaves does not exist in general. Remark that when all the vector fields of $X$ are zero, then $X$ (with the identity map) is a coarse space of leaves. Remark also  that, when it exists, the coarse space of leaves is unique, up to a unique isomorphism. Let us raise a question: 


%

\begin{question}
 Let $R$ be a \D-ring. Let us consider the morphism of \D-schemes $t_R : \DSp R \longto (\Sp C(R), \vec{\mathbf{0}})$ induced by the inclusion $C(R)\subset R$. Is $t_R$ a coarse space of leaves for $\DSp R$ ?
\end{question}
\noindent If $R$ is a constant ring, if $R$ is simple or quasi-simple, then the answer is yes. However, I have some doubt on the positive answer in the general case.

\smallskip

\subsection{Construction of the coarse space of leaves for quasi-simple \D-schemes}

Let $\mathscr{X}$ be a quasi-simple scheme. Let us consider
$\mathcal{O}^{\Delta}_X(X)$
the constant ring of the \D-ring of global sections of $\mathscr{X}$. Let $U$ be a nonempty affine open set of $X$. We have seen that $\mathcal{O}_X(U)$ is quasi-simple. By Proposition \Ref{Cst_faisc_cst_quasi_simple}, the restriction map $\mathcal{O}^{\Delta}_X(X)\longto \mathcal{O}^{\Delta}_X(U)$ is an isomorphism, and by quasi-simplicity, the map $C(\mathcal{O}_X(U)) \longto C (\Frac \mathcal{O}_X(U))$ is also an isomorphism. Hence, $\mathcal{O}^{\Delta}_X(X)$ is a field, seen as a \D-ring by endowing it with  zero derivations. Now, let us construct of morphism 
\[ 
t:\mathscr{X} 
\longto \DSp ( \mathcal{O}_X^{\Delta}(X)).
 \]
 As a continuous function, it is mapping every element of $X$ to the unique element of $\Sp (\mathcal{O}_X^{\Delta}(X))$. Then, the action on scheaves $\mathcal{O}_{\Sp (\mathcal{O}_X^{\Delta}(X))}\longto t_* \mathcal{O}_X$ is simply given by the inclusion map $\mathcal{O}_X^{\Delta}(X)\subset \mathcal{O}_X(X)$, which is a morphism of \D-rings.
We denote 
\[  C(\mathscr{X}):= \DSp ( \mathcal{O}_X^{\Delta}(X))\qquad \aand \qquad t_{\mathscr{X}}: \mathscr{X}\longto C(\mathscr{X})\]
the constructed morphism.

\smallskip

\begin{theorem}\label{Existence_Coarse_Sp_Leaves_Qu_Simple}
Let $\mathscr{X}$ be a quasi-simple \D-scheme. Then, $t_{\mathscr{X}}: \mathscr{X}\longto C(\mathscr{X})$ is coarse space of leaves for $\mathscr{X}$.
\end{theorem}
\smallskip

To prove this theorem, we need to prove a lemma on (unique) factorization and localization. For this result, we set in the category of (\D-)ringed space. In a way, this lemma allows oneself to patch together universal objects. We do not do these patchings directly in the proof of Theorem \Ref{Existence_Coarse_Sp_Leaves_Qu_Simple}, because it would be too much complicated. Hence, we first prove this lemma. We consider  three (\D-)ringed spaces $X$, $Y$ and $Z$, with morphisms $f: X\longto Y$ and $p: X \longto Z$. We want to investigate the existence and the uniqueness of factorizations $g$
\[ 
\xymatrix@R=5mm{
X \ar[r]^{f}\ar[d]_p & Y \\
Z \ar[ur]_g
}
 \]
in terms of (unique) factorizations of ``subsystems'' of $X\longto Y$. We denote
\[ 
E(X, Y):= \left \{  g : Z \longto Y \,\Bigg | \,  
\boite{
\xymatrix@R=5mm{
X \ar[r]^{f}\ar[d]_p & Y \\
Z \ar[ur]_g}
}\text{ commutes}
\right \}
 \]
and $e(X,Y):=\#\, E(X,Y)$.

\begin{lemma}\label{lemme_hyper_technique}
Let $X$, $Y$ and $Z$ be three (\D-)ringed spaces, with $X$ irreducible, and let  $f: X\longto Y$ and $p: X \longto Z$ be morphisms.
\begin{itemize}
\item[(i)] Let $(U_i)_{i\in I}$ be a  basis of nonempty open sets of $X$ such that for all $i\in I$, $e(U_i, Y)=1$. Then, $e(X,Y)=1$.
\item[(ii)] Assume $\# Z=1$. Let $(V_j)_{j\in J}$ be a basis of nonempty sets of $Y$ such that for all $j\in J$ such that $f^{-1}(V_j)\neq \varnothing$, $e(f^{-1}(V_j), V_j)=1$.   Then, $e(X,Y)=1$.
\end{itemize}
\end{lemma}

\pproof{To begin with, let us list some easy properties of $E(-,-)$.
\begin{itemize}\setlength{\itemsep}{0mm}
 \item[\emph{(p1)}] If $U\subset V$ are nonempty open subsets of $X$, then $E(U, Y)\supset E(V, Y)$.
  \item[\emph{(p2)}] If $(W_\ell)_{\ell\in L}$ is a family of nonempty subsets of $X$, then
  \[
  \bigcap_{\ell\in L} E(W_\ell, Y) = E \Big( \bigcup_{\ell \in L}W_{\ell}, Y \Big ).
  \]
\end{itemize}
\smallskip

For \emph{(i)}, to begin with, let us deduce from $(\fa i\in I,\, e(U_i,Y)=1)$ that  all the $E(U_i, Y)$ are equal. Let $i, j\in I$. Let us write the factorizations
\[ 
\boite{
\xymatrix@R=5mm{
U_{i}\ar[r]^{f}\ar[d]_p & Y \\
Z \ar[ur]_{g_1}}
}\quad\aand \quad
\boite{
\xymatrix@R=5mm{
U_{j} \ar[r]^{f}\ar[d]_p & Y \\
Z \ar[ur]_{g_2}}
}.
 \]
Since $X$ is irreducible, $U_{i}\cap U_{j}$ is nonempty and contains $U_{\ell}$ for some $\ell$. One checks then than $g_1$ and $g_2$ belong to $E(U_\ell, Y)$. Hence, $g_1=g_2$ and $E(U_i, Y)=E(U_j, Y)$. Then, one can conclude by \emph{(p2)}.

\smallskip Let us prove \emph{(ii)} now.  In this proof, let us call $j\in J$ \emph{admissible} is  $f^{-1} (V_j)\neq \varnothing$. For all $j$ admissible, let us denote by $g_j$ the unique element of $E(f^{-1}(V_j), V_j)$:
\[
\xymatrix@R=5mm{
f^{-1}(V_j)\ar[r]^{f}\ar[d]_p & V_j \\
Z \ar[ur]_{g_j}}.
\]
Let us denote by $z$ the unique element of $Z$. First, remark that all the $g(z)$ are equal, for all nonempty open set $U$ of $X$ and all $g\in E(U,Y)$. Indeed, if $U$ and $U'$ are two such sets, and if $g\in E(U, Y)$ and $g'\in E(U', Y)$, then $U\cap U'$ is nonempty and contains some $x$. Then, $g(z)$ and $g'(z)$ must be equal to $f(x)$. 
Now, if $j$ is admissible, the inclusion $V_j\subset Y$ gives rise to an injection
\[
E(f^{-1}(V_j), V_j) \longto E(f^{-1}(V_j), Y). 
\]
This map is actually a bijection. Indeed, if $g\in  E(f^{-1}(V_j), Y)$, $g(z)=g_j(z)\in V_j$, and so $g$ can be factorized through $V_j\subset Y$. Hence,
\begin{equation}\label{egalite_E}
 \text{for all $j$ admissible} \qquad e(f^{-1}(V_j), Y)=1.
\end{equation}
Finally, let us prove that $E(f^{-1}(V_i), Y)=E(f^{-1}(V_j), Y)$ whenever $i,j$ are admissible. By \emph{(p2)}, it will conclude this proof. Since, $f^{-1}(V_i)\cap f^{-1}(V_j)$ is nonempty, one can find $\ell$ admissible such that $V_\ell \subset V_i\cap V_j$. By \emph{(p1)}, one has $E(f^{-1}(V_i), Y)\subset E(f^{-1}(V_\ell), Y)$. But, by (\ref{egalite_E}), the cardinal of these both sets is $1$. Hence, $E(f^{-1}(V_i), Y)= E(f^{-1}(V_\ell), Y)$ and so $E(f^{-1}(V_i), Y)= E(f^{-1}(V_j), Y)$. }

\bigskip

\pproofbis{of Theorem \Ref{Existence_Coarse_Sp_Leaves_Qu_Simple}}{
Let $\mathscr{X}$ be a quasi-simple \D-scheme. Let $Y$ be a scheme and let $f : \mathscr{X}\longto (Y, \vec{\mathbf{0}})$ be a morphism. We will also denote $ (Y, \vec{\mathbf{0}})$ by $Y$. With $t_{\mathscr{X}} : \mathscr{X}\longto C(\mathscr{X})$, we are in the situation
\[
\xymatrix@R=5mm{
\mathscr{X}\ar[r]^{f}\ar[d]_{t_{\mathscr{X}}} & Y \\
C(\mathscr{X})
}.
\]
described above. We want to prove that $e(X, Y)=1$.
\smallskip

First, let us assume $Y$ affine, $Y:= \Sp A$. Let $\mathscr{U}$ be a nonempty open affine set of $\mathscr{X}$, endowed with the induced vector fields. Let us write $\mathscr{U}:= \DSp B$. By Proposition \Ref{Cst_faisc_cst_quasi_simple}, one can write $C(\mathscr{X})=\Sp C(B)$. We are now in the situation
\[
\xymatrix@R=5mm{
\DSp B \ar[r]^{f}\ar[d]_{p} & \Sp A \\
\Sp C(B)
}.
\]
Let us prove that $e(\DSp B, \Sp A)=1$. It reduces to prove that given a morphism $\varphi: A \longto B$ of \D-rings (with $A$ endowed with the zero derivations), $\varphi$ factors uniquely through the inclusion $C(B) \subset B$, which is obvious. Hence, $e(\DSp B, \Sp A)=1$ and, by Lemma \Ref{lemme_hyper_technique}.\emph{(i)}, we know now that whenever $Y$ is affine, one has $e(X, Y)=1$.
\smallskip

Now, the general case is easy. We will apply Lemma \Ref{lemme_hyper_technique}.\emph{(ii)} with the basis $(V_i)_i$ of nonempty open affine sets of $Y$. Given such a set $V_i$, if $f^{-1} (V_i)\neq \varnothing $, by Lemma \Ref{Carac_Sch_Quasi_Simple}.\emph{(ii)}, it is quasi-simple, and $C(\mathscr{X})=C(f^{-1}( V_i))$. Thus, one can apply the previous case: $e(f^{-1}( V_i), V_i)=1$. One then concludes with  Lemma \Ref{lemme_hyper_technique}.\emph{(ii)}. 
}

\bigskip

\section{The affine statements}
\label{Affine_Statements}

In this section, we state and prove the generalization of Lemma S-vdP, as a statement of commutative differential algebra. Actually, our geometrization of Lemma S-vdP will be based on two affine statements, Proposition A and Proposition B. Let us recall:
\begin{lemmaSvdP}
Let $k$ be a differential field, with field of constants  $C$. 
The $k$-algebra $k[X_{ij}, \frac{1}{\det X_{ij}}]$ is considered as a differential ring with the derivation defined by ${X_{ij}}'=0$. Then, the maps
$$I \longmapsto \left \langle I  \right \rangle\qquad \aand \qquad J \longmapsto J \cap C[X_{ij}, \frac
{1}{\det X_{ij}}]$$
are inverse bijections between the set of ideals $I$ of $C[X_{ij}, \frac{1}{\det X_{ij}}]$ and the set of  differentials ideals $J$ of $k[X_{ij}, \frac{1}{\det X_{ij}}]$.
\end{lemmaSvdP}

\subsection{The first statement}

If $R$ is a \D-ring, and if $I$ is a subset of $R$, we will denote by $\langle I \rangle$ the \D-ideal of $R$ generated by $I$.

\begin{PropositionA}\label{Gene_Aff_1}
 Let $K$ be a simple \D-ring and $C$ its field of constants. Let $A$ be a $C$-algebra, endowed with the zero derivations. Denote by $i : A \longto A\otimes_C K$ the canonical morphism. Then, the maps
$$I \longmapsto \left \langle i(I)  \right \rangle\qquad \aand \qquad J \longmapsto i^{-1}(J) $$
are inverse bijections between the set of ideals $I$ of $A$ and the set of  \D -ideals $J$ of $A\otimes_C K$.
\end{PropositionA}

\bigskip

Before proving this proposition, let us make several remarks:
\begin{itemize}
\item[\emph{(r1)}] The family $(\widetilde{\partial})_{\partial\in \Delta}$ of derivations of $A\otimes_C K$ is defined by
\[
\fa a\in A,\,  \fa \lambda \in K,\qquad\widetilde{\partial}(a\otimes\lambda):= a \cdot \partial(\lambda)
\]
for all $\partial\in \Delta$. Similarly, given any $\theta\in \Theta(K)$, if one denotes by  $\widetilde{\theta}$ the same operator viewed in $\Theta(A\otimes_C K)$, the one has  $\widetilde{\theta}(a\otimes\lambda)= a \cdot \theta(\lambda)$ for all $a\in A$ and $\lambda \in K$.

 \item[\emph{(r2)}] Given an ideal $I$ of $A$, then the ideal generated by $i(I)$ in $A\otimes_C K$ and the \D-ideal generated by $i(I)$ are equal. Indeed, any $x\in (i(I))$ can be written
 \[
 x=\sum_{j=1} a_j\otimes \lambda_j,
 \]
 where $a_j\in I$ and $\lambda_j\in K$. Then, one computes for all $\widetilde{\theta}\in \Theta(A\otimes K)$,
 \[
 \widetilde{\theta}(x)=\sum_{j=1} a_j\otimes \theta(\lambda_j),
 \]
 which belongs to $(i(I))$.

 \item[\emph{(r3)}] The following two inclusions are easy. First, given an ideal $I$ of $A$, then $I \subset i^{-1}\langle i(I) \rangle$.  Second, given any ideal (not necessarily a \D-ideal) $J$ of $A\otimes_C K$, then $\left ( i(i^{-1} J)\right )\subset J$.

 \item[\emph{(r4)}] Given a \D-ring $R$, we define $R[\boldsymbol{\partial}]$ to be the free $R$-module generated by $\Theta(R)$. It is also a non-commutative $R$-algebra, the product being defined by $\partial \cdot x=\partial (x) +x \cdot\partial$, for all $\partial \in \Delta$. There is a left action of $R[\boldsymbol{\partial}]$ on $R$ defined by
 \[
 \Big (\sum_{\theta\in\Theta(R)} a_\theta \cdot \theta \Big ) \bullet x :=\sum_{\theta\in\Theta(R)} a_\theta \cdot \theta(x).
 \]
 A subset $I$ of $R$ is a \D-ideal iff it is stable under the action of $R[\boldsymbol{\partial}]$.
  Given $x\in R$, the map
 \[
 \text{ev}_x : \fonctionb{ R[\boldsymbol{\partial}]}{R}{L}{L \bullet x}
 \]
 is $R$-linear and is a $R[\boldsymbol{\partial}]$-equivariant map.
We define
 \[
 \text{EqDiffLin}_R(x):= \Ker (\text{ev}_x)=\{ L \in R[\boldsymbol{\partial}] \tq L \bullet x=0 \}.
 \]
 It is the left ideal of $R[\boldsymbol{\partial}]$ of linear differential equations satisfied by $x$. Dually, the \D-ideal generated by $x$ equals $\Im (\text{ev}_x)$.
 
   \item[\emph{(r5)}] 
  If $L=\sum_\theta a_\theta \cdot \theta \in K[\boldsymbol{\partial}]$, we will denote by $\widetilde{L}\in (A\otimes_C K)[\boldsymbol{\partial}]$ the differential operator
 \[
 \widetilde{L}:=\sum_\theta a_\theta \cdot \widetilde{\theta}.
 \]
Then, for $a\in A$ and $\lambda\in K$, one has
$\widetilde{L}\bullet(a \otimes \lambda)=a \otimes (L\bullet \lambda)$.
 
  \item[\emph{(r6)}] Any $x\in A \otimes_C K$ can be written
  $$
  x=\sum_{i=1}^n a_i \otimes \lambda_i,
  $$ 
  with $n\in \N_{\geq 0}$, $(a_i)_i \in A^n$ and $(\lambda_i)_i \in K^n$. We call \emph{length of $x$}, and we denote it by $\ell(x)$ the least $n\geq 0$ such that $x$ can be written as a sum with $n$ terms. One has $\ell(x)=0$ iff $x=0$ and $\ell(x)=1$ iff $x$ is a nonzero pure tensor.
\end{itemize}

\bigskip

\pproofbis{of Proposition A}{Because of \emph{(r3)}, we just have two inclusions to prove. First, given $I$ an ideal of $A$, let us prove that $i^{-1}\langle i(I) \rangle \subset I$. As $C$ is a field, let $(e_j)_{j\in J}$ be a $C$-basis of $K$, with $e_{j_0}=1$. If we denote $\widetilde{e}_j:=1\otimes e_j \in A\otimes_C K$, one knows (Proposition 4.1 of \Cite{algebra}, page 623) that $(\widetilde{e}_j)_{j\in J}
$ is a $A$-basis of $A\otimes_C K$. Let $a\in i^{-1}\langle i(I) \rangle $: there exist $a_1, \ldots, a_n\in I$ and $\lambda_1, \ldots, \lambda_n\in K$ such that
\[
a \otimes 1 =\sum_{i=1}^n a_i\otimes \lambda_i.
\]
Decomposing the $\lambda_i$'s in the $C$-base $(e_j)_j$ one thus finds  $\widetilde{a}_1, \ldots, \widetilde{a}_m\in I$ and an injective map $j : \{1, \cdots, m \}\longto J$ such that
\[
a \otimes 1 =\sum_{i=1}^n \widetilde{a}_i\otimes e_{j(i)} \qquad \qquad \ie, \quad a \cdot \widetilde{e}_{j_0}-\sum_{i=1}^m \widetilde{a}_i\cdot \widetilde{e}_{j(i)}=0.
\]
Since $(\widetilde{e}_j)_j$ is a $A$-basis, $a$ must be equal to one of the $\widetilde{a}_i$ and so $a\in I$.
\bigskip

Now, given $J$ a \D-ideal of $A\otimes_C K$, let us prove that $J\subset \langle i(i^{-1} J) \rangle$. More precisely, let us prove by induction on $\ell(x)$ that
\begin{equation}\label{hyp_induction}
 x\in J \quad \impl \quad  x\in\langle i(i^{-1} J) \rangle.
\end{equation}
If $\ell(x)=0$, this means that $x=0$, and (\ref{hyp_induction}) is clear. If $\ell(x)=1$, then one can write $x=a \otimes \lambda$, with $a\in A$ and $\lambda\in K$ nonzero. Since $\lambda\neq 0$ and $K$ is simple, $\langle \lambda \rangle=K$. Hence, there exists $L\in R[\boldsymbol{\partial}]$ such that $L \bullet \lambda =1$. Then, $\widetilde{L}\bullet x = a\otimes 1 \in J$, since $J$ is a \D-ideal. Thus, $a\in i^{-1} (J)$, and $x\in \langle i(i^{-1} J) \rangle$. Now, assume (\ref{hyp_induction}) true whenever $\ell(x)\leq n$ and let $x\in J$ with $\ell(x)=n+1$. Let $(a_i)_i\in A^{n+1}$ and  $(\lambda_i)_i\in K^{n+1}$ such that
\[
x=\sum_{i=1}^{n+1} a_i \otimes \lambda_i. 
\]
First, remark that one can suppose that
\begin{equation}\label{hyp_EqDiffLin}
 \fa (i,j)\in \{1, \ldots, n+1\}^2, \qquad  \text{EqDiffLin}_K(\lambda_i)= \text{EqDiffLin}_K(\lambda_j).
\end{equation}
Indeed, let us suppose that it is not the case, and let $i_1$ and $i_2$ such that
\[
 \text{EqDiffLin}_K(\lambda_{i_2})\subsetneq \text{EqDiffLin}_K(\lambda_{i_1}).
\]
Then, let $L\in K[\boldsymbol{\partial}]$ such that $L\bullet \lambda_{i_1}= \widehat{\lambda_{i_1}}\neq 0$ and $L\bullet \lambda_{i_2}= 0$. Actually, applying another differential operator, one can assume that $L \bullet \lambda_{i_1}=\lambda_{i_1}$ and $L\bullet \lambda_{i_2}=0$. Now, $\widetilde{L}\bullet x \in J$ and $\ell(\widetilde{L}\bullet x )\leq n$, since $L$ kills one of the $\lambda_i$'s. Thus, $\widetilde{L}\bullet x  \in  \langle i(i^{-1} J) \rangle$. But, $x- \widetilde{L}\bullet x$ has also a length $\leq n$ and so $x- \widetilde{L}\bullet x$ and $x$ belong to $\langle i(i^{-1} J) \rangle$, what we wanted.

\smallskip

So, let us suppose (\ref{hyp_EqDiffLin}) holds. Since $\lambda_1\neq 0$, let $L\in R[\boldsymbol{\partial}]$ such that $L \bullet \lambda_1=1$. By (\ref{hyp_EqDiffLin}), one easily gets that all the $L \bullet \lambda_i$'s are constant. Let us denote $c_i:=L \bullet \lambda_i$. Then, one has
\[
\widetilde{L}\bullet x = \sum_{i=1}^{n+1} a_i \otimes (L\bullet \lambda_i) = \sum_{i=1}^{n+1} a_i \otimes c_i =\Big ( \sum_{i=1}^{n+1} a_ic_i \Big ) \otimes 1. 
\]
Hence, $\ell(\widetilde{L}\bullet x )\leq 1$ and so $\widetilde{L}\bullet x \in \langle i(i^{-1} J) \rangle$. But, $x- \lambda_1 \cdot \widetilde{L}\bullet x$ has also length $\leq n$ and so belongs to $\langle i(i^{-1} J) \rangle$. Thus, $x\in \langle i(i^{-1} J) \rangle$.
}

\subsection{The second statement}

\begin{PropositionB}
 Let $K$ be a simple \D-ring, with field of constants $C$. Let $A$ be a $C$-algebra, endowed with the zero derivations. Let $\varphi : K \longto C$ be a morphism of $C$-algebras. Let us denote
 \[
 i : \fonctionb{A}{A\otimes_C K}{a}{a\otimes 1} \qquad \aand \qquad j:\fonctionb{A\otimes_C K}{A}{a\otimes \lambda}{a\varphi(\lambda)}.
 \]
 Then for all \D-ideal $J$ of $A\otimes_C K$ and for all ideal $I$ of $A$, one has
 \[
 i^{-1}(J)=j(J) \qquad \aand \qquad j^{-1}(I)_{\#}=(i(I)).
 \]
\end{PropositionB}

\smallskip

Before proving this proposition, let us make two remarks:

\begin{itemize}
 \item[\emph{(r7)}]  Since $j$ is surjective, $j(J)$ is  already an ideal for all ideal $J$ of $A\otimes_C K$. Thus, for all ideal $I$ of $A$, $j(j^{-1} I)=I$.
 \item[\emph{(r8)}] Given any ideal $J$ (not necessarily a \D-ideal) of $A\otimes_C K$, one has $i^{-1} (J) \subset j(J)$. 
\end{itemize}

\smallskip

\pproof{
Let us prove the first assertion. Let $J$ be a \D-ideal of $A\otimes_C K$. By \emph{(r8)}, it suffices to check that $j(J)\subset i^{-1}(J)$. Let $y:=j(x)$ be some element of $j(J)$, for some $x\in J$. By Proposition~A, one has $x\in (i(i^{-1}(J))$. This means that one can write
\[
x=\sum_{i=1}^n a_i \otimes \lambda_i
\]
with the $a_i$'s satisfying $a_i\otimes 1\in J$. Thus,
\[
y=\sum_{i=1}^n a_i \cdot \varphi(\lambda_i).
\]
Since  $a_i\otimes 1\in J$, $a_i\in i^{-1}(J)$ and so $y \in i^{-1}(J)$.
\bigskip

Now, let us turn to the second assertion. Let us recall that, in a \D-ring $R$, given an ideal $\mathfrak{a}$ of $R$, $\mathfrak{a}_\#$ denotes the largest \D-ideal contained in $\mathfrak{a}$. First, check that $(i(I))\subset j^{-1}(I)$. Then, by \emph{(r2)}, $(i(I))$ is a \D-ideal, and so one has  $(i(I))\subset j^{-1}(I)_{\#}$. So, we are asking if $(i(I))$ is the largest \D-ideal contained in $j^{-1} (I)$. 
Hence, let $J$ be a \D-ideal such that
\[
(i(I))\subset J \subset j^{-1} (I).
\]
Applying $i^{-1}$ to these inclusions, one gets
\[
i^{-1} ((i(I))) \subset i^{-1} (J) \subset i^{-1} (j^{-1} (I)).
\]
By Proposition A, $i^{-1} ((i(I)))=I$ and by \emph{(r8)}, $i^{-1} (j^{-1} (I)) \subset j(j^{-1}(I))$. By \emph{(r7)}, $j(j^{-1}(I))=I$. So, one has
\[
I \subset i^{-1} (J) \subset I \qquad \text{and so} \qquad i^{-1} (J)=I.
\]
Applying $(i(-))$ to these identity, one gets
\[
(i(I))= (i(i^{-1} (J))) =J
\]
again by Proposition A. Hence, $(i(I))=j^{-1}(I)_{\#}$. }

\bigskip

\begin{corollary}
Let $K$ be a simple $\Q$-\D-algebra, with field of constants $C$. We assume there exists a morphism $\varphi : K \longto C$. Let $A$ be a $C$-algebra, endowed with the zero derivations. We consider $A \otimes_C K$ with the canonical morphism $i: A\longto A\otimes_C K$. Then, for all $\mathfrak{p}\in \Sp A$, the ideal $(i(\mathfrak{p}))$ is a prime \D-ideal.
\end{corollary}

\pproof{
Use Proposition B and Proposition \Ref{primality_ideals}.\emph{(iv)}.
}
\bigskip

Can one drop the hypothesis that $K$ is defined over $\Q$, and that $\Sp K$ has a $C$-point ? It may be possible, by using similar techniques as in the proof of Proposition A. More explicitely, we ask the question:

\begin{question}
Let $K$ be a simple \D-ring, with $C$ as field of constants. Let $A$ be a $C$-algebra and let $i:A\longto A\otimes_C K$ the canonical morphism. Does one have
\[
\mathfrak{p}\in \Sp A \ \impl \ (i(\mathfrak{p}))\in \Sp (A\otimes_C K) ?
\]
\end{question}

\bigskip

\section{Trajectories of \D-schemes without self-dynamics over a simple base}

We are now able to prove our geometrization of Lemma S-vdP. To begin with, let us fix some notations that  will  be  used throughout this section.

\begin{itemize}
 \item $\mathscr{S}$ is simple \D-scheme. The family of vector fields attached to $\mathscr{S}$ is denoted by $\fchv{V}_S$.

 \item $t_{\mathscr{S}} : \mathscr{S} \longto C(\mathscr{S})$ is a coarse space of leaves for $\mathscr{S}$.

  \item$\mathscr{X}\longto \mathscr{S}$ is a \D-scheme above $\mathscr{S}$, without self-dynamics (see below for the definition). One denotes by $X_0$ a model of $\mathscr{X}$. The vector fields attached to $\mathscr{S}$ are denoted by $\fchv{V}$.

  \item$x : C(\mathscr{S}) \longto \mathscr{S}$ is a $C(\mathscr{S})$-point of $\mathscr{S}$, compatible with the $C(\mathscr{S})$-structure of $\mathscr{S}$, that is 
  $\xymatrix@1{C(\mathscr{S})\ar[r]^-{x}
  &\mathscr{S}\ar[r]^-{t_{\mathscr{S}}}&C(\mathscr{S})}$ is the identity.

  \item $X_x$ will denote the fiber of $\mathscr{X}$ above $X$. The (closed) immersion of $X_x$  into $X$ will be denoted by $i_{X,x} : X_x \longto X$. 
  
  \item If $T$ is an open set of $S$, we denote by $\mathscr{X}_T$ the preimage of $T$ by $\mathscr{X} \longto \mathscr{S}$. Equivalently, the square 
  \[
  \xymatrix@R=6mm@C=7mm
{
\mathscr{X}_T  \immouv[r] \ar[d] & \mathscr{X} \ar[d] \\
 T \immouv[r]  &\mathscr{S}
}
  \]
  is cartesian.

\item If $U_0$ is an open set of $X_0$, $\mathscr{U}$ will denote the open subscheme of $\mathscr{X}$ defined by the cartesian square
\[
\boite{
\xymatrix@=7mm{
\mathscr{U}\immouv[d] \ar[r] & U_0\immouv[d] \\
\mathscr{X} \ar[r] & X_0
}
}.
\]

    \item When affine, we will denote $\mathscr{S}=\DSp K$, $C(\mathscr{S})=\Sp C$, where $K$ is a simple \D-ring with field of constants $C$, $X_0=\Sp A$, $\mathscr{X}=\DSp (A\otimes_C K)$, with $A$ any $C$-algebra.
\end{itemize}

\smallskip

\noindent In this section, we will use more or less implicitely these three facts:
\begin{itemize}
 \item[\emph{(f1)}] In the category of schemes, open immersions are stable under pull-back. 
  \item[\emph{(f2)}] In any category, consider objects and arrows
  \[
  \xymatrix@=4mm{
  && X\ar[d] \\
  \widetilde{S}\ar[r] & S' \ar[r]&S
  }
  \]
  and assume you want to compute a fibered product $\widetilde{X}:=X\times_S \widetilde{S}$. Then, it is sufficient to compute first a fibered product $X':=X\times_{S} S'$ and second $\widetilde{X'}:=X'\times_{S'}\widetilde{S}$.  In other words, $X\times_S \widetilde{S}=(X\times_S S' )\times_{S'} \widetilde{S}$.
    \item[\emph{(f3)}] Let 
      \[
  \xymatrix@=4mm{
  &Y\ar[d] \\
  W \ar[r]&Z
  }
  \]
 be a diagram of schemes and let $(U_i)_{i\in I}$ an open covering of $Y$. Then, $(U_i\times_Z W)_{i\in I}$ is an open covering of $W$.  This follows from \Cite{EGA1} Proposition (4.3.1) and Proposition (3.5.2). A nice way to understand this fact is to say that the Zarisiki coverings define a Grothendieck pretopology. \end{itemize}

\subsection{\D-schemes without self-dynamics}

In \Cite{buium}, Buium defines \emph{a $K$ split \D-variety} to be a \D-variety that one can descend to $C(K)$. With the help of coarse space of leaves, we are able to give a similar definition for non-affine schemes.

\begin{definition}
Let $\mathscr{S}$ be a \D-scheme that has a coarse space of leaves $t_{\mathscr{S}} : \mathscr{S}\longto C(\mathscr{S})$. Let $\mathscr{X}$ be a $\mathscr{S}$-\D-scheme. We say that \emph{$\mathscr{X}$ has no self-dynamics} if there exists $X_0$ a $C(\mathscr{S})$-scheme and a morphism $\mathscr{X}\longto (X_0,\vec{\mathbf{0}})$ making the square
\[ 
\xymatrix@=5mm{
\ar[d]\mathscr{X}\ar[r]& (X_0, \vec{\mathbf{0}})\ar[d]\\
\mathscr{S}\ar[r] & C(\mathscr{S})
}
 \]
 cartesian.
\end{definition}

\subsection{Intersection of a leaf with the fiber}

\begin{proposition}\label{prop_intersection}
Assume $\mathscr{S}$ to be of characteristic zero.
 Let $\eta\in {X}^{\fchv{V}}$ a leaf of $\mathscr{X}$. Then, $\eta$ intersects $X_x$ in a unique point. That is, there exists a unique $\eta_{\mid x}\in X_x$ such that
 \[
 \fa y\in X_x, \qquad \eta \leadsto i_{X,x}(y) \ssi \eta_{\mid x} \leadsto y.
 \]
\end{proposition}

\begin{lemma}\label{Lemma_appartenance_feuille}
  Let $\eta\in \mathscr{X}^{\fchv{V}}$ a leaf of $\mathscr{X}$. Then, for every nonempty open set $T$ of $ \mathscr{S}$, $\eta \in \mathscr{X}_T$.
\end{lemma}

\pproof{First, let us prove that one can assume that $\mathscr{X}$ is affine. Let $\eta$ be a leaf of $\mathscr{X}$. By \emph{(f3)}, there exists an open affine set $U_0$ of $X_0$ such that $\eta\in \mathscr{U}$. If we prove that for all nonempty set $T$ of $S$, $\eta\in \mathscr{U}_T$, then since $\mathscr{U}_T$ is an open subset of $\mathscr{X}_T$, we will have proved the lemma. So, in the following, we assume $\mathscr{X}$ affine.

\smallskip

Then, let us prove this lemma when $\mathscr{S}$ is affine. Let $K$ be a simple \D-ring, $C$ its field of constants, $A$ a $C$-algebra. We can assume furthermore that $T$ is a distinguished open set $D(f)$ for some $f\in K$. Let $\mathfrak{p}$ be a prime \D-ideal of $A\otimes_C K$. Since in the diagram
\[
  \xymatrix@R=6mm@C=7mm
{
\mathscr{X}_T  \immouv[r] \ar[d] & \mathscr{X} \ar[d] \ar[r] &X_0 \ar[d]\\
 T \immouv[r]  &\mathscr{S} \ar[r] & C(\mathscr{S})
}
\]
all the square are cartesian, the arrow $\mathscr{X}_T \longto \mathscr{X}$ can actually been described as the spectrum of
\[
i_f: A\otimes_C K \longto A \otimes_C K_f,
\]
and one wants to prove that there exists $\mathfrak{q}\in \Sp (A \otimes_C K_f)$ such that $\mathfrak{p}={i_f}^{-1}( \mathfrak{q})$. Now, if $x\in  A\otimes_C K$ and $f^n x\in \mathfrak{p}$, then $x\in  \mathfrak{p}$. Indeed, if $x\notin \mathfrak{p}$, then $1\otimes f^n\in \mathfrak{p}$ and as $f^n \neq 0$, let $L\in K[\boldsymbol{\partial}]$ such that $L \bullet (f^n)=1$. Since $\mathfrak{p}$ is a \D-ideal, $\widetilde{L}\bullet (1 \otimes f^n)=1\in \mathfrak{p}$, which is absurd. Let us set $\mathfrak{q}:= i_f(\mathfrak{p})$ and check that $\mathfrak{q}$ is prime. Let $x/f^n$ and $y/f^m$ be elements of $A\otimes_C K_f$ such that $xy/f^{n+m}\in \mathfrak{q}$. This means there exist $p\in \N$ such that $f^pxy\in \mathfrak{p}$. By what preceeds, it implies $xy\in \mathfrak{p}$ and so $\mathfrak{q}$ is prime.  Similarly, one proves that $(i_f)^{-1}(\mathfrak{q})=\mathfrak{p}$.
\smallskip

Now, let us turn to the general case. Let $T$ be a nonempty open set of $\mathscr{S}$. 
By \emph{(f3)}, one knows that the $\mathscr{X}_V$'s for $V$ nonempty open sets of $S$ form a open covering of $\mathscr{X}$. Let $V_0$ be a nonempty affine set such that $\eta\in  \mathscr{X}_{V_0}$. Clearly, it is still a leaf. By the affine case, and since $T \cap V_0\neq \varnothing$ (for $S$ is irreducible), $\eta$ belongs to $ \mathscr{X}_{(T\cap V_0)}$ and so $\eta\in \mathscr{X}_T$.
}

\bigskip

\pproofbis{of Proposition \Ref{prop_intersection}}{This proof is organized as follows:
\emph{a)} Proof of the unicity.
  \emph{b)} One can assume the base $S$ affine.
  \emph{c1)}  An intermediate result on the injectivity of $\mathscr{X}^\fchv{V}\longto X_0$.
  \emph{c2)} An intermediate result on ``split'' open sets that contains $\eta$.
   \emph{d)} One can assume the total space $X$ affine.
       \emph{e)} Proof of the proposition when $X$ and $S$ are affine.
\bigskip

\emph{a)} Unicity is easy. Assume one has two elements $\eta_{\mid x, 1}$ and $\eta_{\mid x, 2}$ satisfying the conclusions of the proposition. Since $\eta_{\mid x, 1}\leadsto \eta_{\mid x, 1}$, one has $\eta \leadsto i_{X,x}(\eta_{\mid x, 1})$ and thus $\eta_{\mid x, 2}\leadsto \eta_{\mid x, 1}$. Similarly, $\eta_{\mid x, 1}\leadsto \eta_{\mid x, 2}$ and so $\eta_{\mid x, 1}= \eta_{\mid x, 1}$.
\smallskip

\emph{b)} Second, one can assume $\mathscr{S}$ affine. Indeed, if $U$ is an affine open set of $S$ containing the image of $x$, then one has
\[
\xymatrix@R=5mm@C=7mm
{
X_x \ar[r]_-{i_1}\ar[d]  & \mathscr{X}_U  \immouv[r] \ar@{}[r]_-{i_2} \ar[d] & \mathscr{X} \ar[d] \\
C( \mathscr{S}) \ar[r] &  U \immouv[r]  &\mathscr{S}
}
\]
in which every square is cartesian. By Lemma \Ref{Lemma_appartenance_feuille}, one knows that $\eta\in \mathscr{X}_U$. This means that there exist $\eta_U\in \mathscr{X}_U$ such that $i_2(\eta_U)=\eta$. Assume the proposition true when $S$ is affine. Then, let $\eta_{\mid x}\in X_x$ such that $\eta_{\mid x} \leadsto y \ssi \eta_U \leadsto i_1(y)$ for all $y\in X_x$. Now, since $i_2$ is an open immersion, $\eta_U \leadsto i_1(y) \ssi i_2(\eta_U) \leadsto i_2(i_1(y))$, that is $\eta_{\mid x} \leadsto y \ssi \eta \leadsto i_{X,x}(y)$ for all $y\in X_x$. Thus, the general case follows from the case when $S$ is affine.
\smallskip

\emph{c1)} We assume $S$ affine. Let us denote, in what follows, $p: \mathscr{X}\longto X_0$ the structure map.  Let us prove that $p_{\mid \mathscr{X}^{\fchv{V}}} : \mathscr{X}^{\fchv{V}} \longto X_0$ is injective. Let $\eta_1$ and $\eta_2$ be two leaves such that $p(\eta_1)=p(\eta_2)$. Let $U_0$ be an open affine of $X_0$ that contains $p(\eta_1)$. Then, $\eta_1$ and $\eta_2$ belong to $\mathscr{S} \times_{C(\mathscr{S})} U_0 $, which is affine. Thus, one can assume $X$ affine. What we want to prove is that given two prime \D-ideals $\mathfrak{p}_1$ and $\mathfrak{p}_2$ of $A\otimes_C K$, then $i^{-1} (\mathfrak{p}_1) = i^{-1}( \mathfrak{p}_2)$ implies $\mathfrak{p}_1 =\mathfrak{p}_2$. This is clear by Proposition A.

\smallskip

\emph{c2)} We assume $S$ affine.  Let $U_0$ be open set of $X_0$. We denote $\mathscr{U}:= \mathscr{S} \times_{C(\mathscr{S})} U_0 $. By \emph{(f1)} and \emph{(f2)}, the diagram
\begin{equation}
\boite{
 \xymatrix@R=5mm@C=7mm
{
\mathscr{U}\immouv[r]\ar[d] & \mathscr{X_0} \ar[d] \\
  U_0 \immouv[r]  &X
}
}\label{un_diagrame_comm}
\end{equation}
is a cartesian square.
Recall that $\eta$ is a leaf of $\mathscr{X}$. We prove that $\eta\in \mathscr{U}$ iff $p(\eta)\in U_0$. The direction $\impl$ is clear by commutativity of (\ref{un_diagrame_comm}). Now, assume $p(\eta)\in U_0$. We know that the map $\mathscr{U} \longto U_0$ is surjective, as it is a base change  of $\mathscr{S}\longto C(\mathscr{S})$ and as surjectivity is stable under base change. So, let $x\in \mathscr{U}$ such that $p(z)=p(\eta)$. Since $\mathscr{U} \longto  (U_0, \vec{\mathbf{0}})$ is a morphism of \D-schemes, we know that $z$ and $\Traj{V}(z)$ have the same image\footnote{Note that this is the only place where we use that $\mathscr{S}$ has characteristic zero.}. Furthermore, by stability of open sets under generization,  $\Traj{V}(z)\in \mathscr{U}$. So, $p(\Traj{V}(z))=p(\eta)$. One concludes with \emph{c1)} that $\eta\in \mathscr{U}$.
\smallskip

\emph{d)} Now, let us prove that one can assume $X_0$ affine. We take $S$ affine and assume that the result holds when $X_0$ is affine. By \emph{(f3)}, the schemes $\mathscr{U}:=\mathscr{S} \times_{C(\mathscr{S})} U_0$ form a open covering of $X$ for $U_0$ nonempty affine open sets of $X_0$. So, let $U_0$ be such a set with $\eta\in \mathscr{U}$. We are in the following situation
\[ 
\boite{
\xymatrix{
X_x \ar[rdd] \ar[r]^{i_{X,x}} & \mathscr{X} \ar[r]^{p}\ar[rdd] |!{[d];[dr]}\hole& X_0\ar[rdd]|!{[d];[dr]}\hole\\
&U_x\immouv[ul]\ar@{}[ul]_<<<<<{j_U} \ar[d] \ar[r]_{i_{U,x} \ \ \ \ }& \mathscr{U}\ar[d] \ar[r]_{p_U \ \ \ \ } \immouv[ul] & U_0 \ar[d] \immouv[ul]  \\
&C(\mathscr{S})\ar[r]_x& \mathscr{S}\ar[r]_{t_{\mathscr{S}}}&C(\mathscr{S})
}
},
 \]
 where all the squares are cartesian. Remark that since the composition of the two arrows on the base line equals $\Id_{C(\mathscr{S})}$, one can assume that $X_x=X_0$ and $U_x=U_0$, with  $i_{X,x}\circ p= \Id_{X_0}$  and $i_{U,x}\circ {p}_U= \Id_{U_0}$.  Now, since we assume the affine case (base and total space both affine) done, there exists $\eta_{\mid U,x}\in U_x$ satisfying the conclusion of the proposition. Let $\mathbb{U}$ be the set of all $U_0$, affine open subsets of $X_{0}$, such that $\eta\in \mathscr{U}$. Let us prove 
\[ 
\fa U_0, V_0 \in \mathbb{U}, \qquad j_U(\eta_{\mid U,x}) =  j_V(\eta_{\mid V,x}).
 \]
 It is easy, indeed, since $U_0\cap V_0\neq \varnothing$, one can find $W_0\in \mathbb{U}$ such that $W_0\subset U_0\cap V_0$. By unicity of $\eta_{W,x}$, one can check that both $\eta_{\mid U,x}$ and $\eta_{\mid V,x}$ equal $\eta_{W,x}$.
 Hence, let us denote by $\eta_{\mid x}$ this unique image of all the $j(\eta_{\mid U,x})$ and let us prove that it is the point we are looking for. 

\smallskip

To begin with, let us prove that
\begin{equation}\label{petite_relation}
p(\eta)\leadsto p(i_{X,x}(\eta_{\mid x})).
\end{equation} 
 It is easy. Let $U_0\in \mathbb{U}$. Then, one has $\eta\leadsto i_{U,x}(\eta_{\mid U,x})$ in $\mathscr{U}$, but also in $\mathscr{X}$. But, $ i_{U,x}(\eta_{\mid U,x})= i_{X,x}(\eta_{\mid x})$, and so, applying $p$, one gets the required relation.

\smallskip
 
 Let $y\in X_x$ such that $\eta\leadsto i_{X,x}(y)$. We want to prove that $\eta_{\mid x}\leadsto y$. First, one has $p(\eta)\leadsto y$. Let $U_0$ be an affine open set such that $y\in U_0$. Hence, $p(\eta)\in U_0$ and so, by \emph{c2)}, $U_0\in \mathbb{U}$. So, the relation $\eta\leadsto i_{X,x}(y)$ can be read in $\mathscr{U}$ as $\eta\leadsto i_{U,x}(y)$ and by definition of $\eta_{U,x}$, one has $\eta_{\mid U,x}\leadsto y$. Applying $j_U$ to this relation, one obtains $\eta_{\mid x}\leadsto y$. Conversely, if $\eta_{\mid x}\leadsto y$ in $X_x$, let us consider an open affine set $U_0$ containing $y$. Then, one has
 \begin{align*}
p(i_{X,x}(\eta_{\mid x})) = \eta_{\mid x} \leadsto p(y),
\end{align*}
 and so, by (\ref{petite_relation}), one has $p(\eta)\leadsto p(y)$. So, $p(\eta)\in U_0$, and by \emph{c2)}, $\eta\in \mathscr{U}$ and so $U_0\in \mathbb{U}$. One easliy checks then than $\eta \leadsto i_{X,x}(y)$.
 
\bigskip

\emph{e)} To end this proof, let us show that the result holds when both the base and the total space are affine. The situation is this case is the following
\[ 
\boite{
\xymatrix{
A \ar@{-}[d]& A \otimes_C K \ar[l]-_{j} \ar@{-}[d]& A \ar[l]_-{i} \ar@{-}[d]\\
C & K \ar[l]& C\ar[l]
}
}.
 \]
 Let $\mathfrak{p}$ be a prime \D-ideal of $A\otimes_C K$. The intersection point we are looking for is  $j(\mathfrak{p})$. It is a prime ideal by Proposition B. The property we want to check is
\[ 
\fa \mathfrak{q}\in \Sp A, \qquad j(\mathfrak{p})\subset \mathfrak{q} \ \ssi \ \mathfrak{p}\subset j^{-1}(\mathfrak{q}).
 \]
This is straightforward. }

\subsection{Trajectory of \D-schemes without self-dynacmics over a simple base}

We can now state and prove our result:

\begin{theorem}\label{main_thm}
Let $\mathscr{S}$ be a simple \D-scheme  and $\mathscr{X}$ be a $\mathscr{S}$-scheme without self-dynamics. Let $\mathscr{S} \longto C(\mathscr{S})$ be a coarse space of leaves of $\mathscr{S}$. Let $x$ be a $C(\mathscr{S})$-point of $\mathscr{X}$. Then, the maps
$$
\fonctionb{X_x}{{X}^{\fchv{V}}}{y}{\Traj{V}(y)}
\qquad \aand \qquad
\fonctionb{X^\fchv{V}}{X_x}{\eta}{\eta_{\mid x}}
$$ 
are inverse bijections between the fibre $X_x$ and the set $X^\fchv{V}$ of leaves of $X$ under the action of $\fchv{V}$.
\end{theorem}

\bigskip

\pproof{First, let us fix $y\in X_x$, and let us prove that $(\Traj{V}(y))_{\mid x}=y$. More precisely, by $\Traj{V}(y)$, we mean $\Traj{V}(i_{X,x}(y))$. Let us check that $y$ satisfies the property of the intersection $\eta_{\mid x}$ of a leaf with $X_x$: we want to show that
\[ 
\fa z\in X_x, \qquad y \leadsto z \ssi \Traj{V}(y)\leadsto i_{X,x}(z).
 \]
The direction $\impl$ follows from $\Traj{V}(y) \leadsto y$. Conversely, if  $\Traj{V}(y)\leadsto i_{X,x}(z)$, by composing with $p$, and since one can assume $p\circ i_{X,x}=\Id$, one obtains $p(\Traj{V}(y))\leadsto z$. But, $p(\Traj{V}(i_{X,x}(y)))=p((i_{X,x}(y))=y$ and so $y \leadsto z$.

\smallskip

Now, let us prove that $\Traj{V}(\eta_{\mid x})=\eta$ for all leaf $\eta$. First, one has  $\eta\leadsto i_{X,x}(\eta_{\mid x})$. Thus, by definition of the trajectory, $\eta \leadsto \Traj{V} (i_{X,x}(\eta_{\mid x}))$. Conversely, as in the proof of Proposition \Ref{prop_intersection}, one can assume $\mathscr{S}$ and $X_0$ affine. Then, one want to check that for all differential prime ideal $\mathfrak{p}$ of $A\otimes_C K$, one has $\mathfrak{p}\subset j^{-1}( j(\mathfrak{p}))_{\#}$. By Proposition B, one has $j^{-1}(I)_{\#}=(i(I))$, for all ideals $I$, and so $ j^{-1}( j(\mathfrak{p}))_{\#}=(i(j(\mathfrak{p})))$. By the same proposition, one has $j(J)=i^{-1}(J)$, and so $ j^{-1}( j(\mathfrak{p}))_{\#}=(i(i^{-1}(\mathfrak{p})))$. By Proposition A, this latter equals $\mathfrak{p}$.
 }

\appendix

\section{Some  commutative \D-algebra}\label{appendix}

In this appendix, we present and prove some results of commutative \D-algebra. Many of these results are known for usual ``differential rings'', that is for ordinary \D-rings. Also, some of these results are known in the case of partial \D-rings. For instance, Kolchin only considers such \D-rings in \Cite{DiffAlgGroupsKolch}. Finally, some of these results have also been proved for Hasse-Schmidt derivations by F. Benoist in \Cite{FranckBenoist}. The results are often straightforward generalizations of the partial case or the ordinary case.
\smallskip

The main goal of this appendix is to provide a proof of Theorem \Ref{theo_traj} in the affine case: for every prime ideal $\mathfrak{p}$, the \D-ideal $\mathfrak{p}_\#$ is also prime --- at least for $\Q$-\D-algebras. This is proved by Keigher in \Cite{KeigherDiffRingFromQuasPrimes}. We explain this proof, based on properties of differential operators in paragraph \RefPar{Proof_keigher}. His proof works only for partial \D-rings. We give a new proof (see Proposition~\Ref{primality_ideals}) of this result in paragraph \RefPar{Prim_D_ideals}, valid for all \D-rings.
\smallskip

Another goal is to provide a study of simple \D-rings.

\subsection{Colon ideals and radical ideals in \D-rings}\label{colon_ideals}

Recall that, in a ring $R$, given an ideal $I$ and a set $S$, one denotes
$
(I:S):= \{x\in R \tq xS \subset I\}$. This is an ideal, called a \emph{colon ideal}. We will also denote
$(I:S^{\infty}):=\{ x\in R \tq \xt n>0, \ xS^{n}\subset I\}$, where $S^{n}$ denotes the set of all possible products $s_1\cdots s_n$, with the $s_i$'s in $S$. It is also an ideal. In the case where $S=\{s\}$, we denote it by $(I:s)$ and $(I:s^{\infty})$. The following lemma gathers some easy properties of colon ideals.

\begin{lemma}\label{lemma_colon_ideal}
Let $R$ be a ring. Let $S$ be a subset of $R$. Let $I,J$ be an ideals of $R$.
\begin{itemize}\setlength{\itemsep}{0mm}
\item[(i)] $I\subset (I:S)\subset (I: S^{\infty})$.
\item[(ii)]$(I:S)\cdot S\subset I$ as sets and $(I:J)\cdot J \subset I$ as ideals.
\item[(iii)] If $R$ is a \D-ring and $I$ a \D-ideal, then $(I:S^{\infty})$ is a \D-ideal.
\item[(iv)] If $I$ is a radical ideal:
\begin{itemize}\setlength{\itemsep}{0mm}
\item[a)] $(I:S)=(I:S^{\infty})$.
\item[b)] $(I:S)$ is a radical ideal.
\end{itemize}
\item[(v)] If $R$ is a \D-ring and $I$ a radical \D-ideal, then $(I:S)$ is a radical \D-ideal.
\end{itemize}
\end{lemma}

\pproof{\emph{(i)},  \emph{(ii)} and \emph{(iv)} are staightforward, and \emph{(iii)} follows from Lemma \Ref{lemma_easy}. \emph{(v)} follows from \emph{(iii)} and \emph{(iv)}.
}

\smallskip

\begin{proposition}\label{prop_min_rad_diff}
Let $R$ be a \D-ring and $I,J$ ideals of $R$. 
\begin{itemize}\setlength{\itemsep}{0mm}
\item[(i)] There exists a minimal radical \D-ideal containing $I$. We denote it by $\{ I \}$. For all $n\geq 0$ and all $\phi_1, \ldots, \phi_n\in \{ \sqrt{-}, \langle - \rangle\}$, one has $\phi_n \circ \cdots \phi_1 (I)\subset\{I\}$.
\item[(ii)] If $I\subset J$ then $\{ I \} \subset \{ J\}$.
\item[(iii)] For all $x,y\in R$, for all  $\theta_1, \theta_2 \in \Theta(R)$, \quad $xy\in I \, \impl \, \theta_1(x)\theta_2(y)\in \{I\}$.
\item[(iv)] $\{I\} \{J\} \subset \{IJ\}$.
\item[(v)] $\{ \{I\} \{J\} \}= \{IJ\}$.
\end{itemize}
\end{proposition}

\pproof{The proof of Lemma 1.6 in \Cite{kaplansky} is still valid for  \D-rings.
For \emph{(i)}, remark that \D-ideals and radical ideals are stable under intersection. For \emph{(iii)}, proceed by induction. The induction step follows from $(x'y)^{2}=(x'y)(xy)'-(xy)(x'y')$.\smallskip

Now, let us prove $\emph{(iv)}$. Let $x\in I$. Then, by \emph{(v)} of Lemma \Ref{lemma_colon_ideal}, the ideal $(\{xJ\} : x)$ is a radical \D-ideal. Furthermore, one checks that $J\subset (\{xJ\} : x)$. Hence, $\{ J\} \subset  (\{xJ\} : x)$. So, by  \emph{(ii)} of Lemma \Ref{lemma_colon_ideal}, $x \{ J\} \subset \{xJ\}$. Now, since $xJ \subset IJ$, one has $x \{ J\} \subset \{IJ\}$. This means that $I\subset (\{IJ\}:\{J\})$. But, $(\{IJ\}:\{J\})$ is a radical \D-ideal. So, $\{I\}\subset (\{IJ\}:\{J\})$. By  \emph{(ii)} of Lemma \Ref{lemma_colon_ideal}, one gets $\{I\} \{J\} \subset \{IJ\}$.

\smallskip The proof of \emph{(v)} is straightforward.
}

\bigskip
\subsection{The monoid of differential operators of a \D-ring in the partial case}\label{Proof_keigher}

 As explained in Subsection \RefPar{D_Rings}, given a partial \D-ring, one has $\Theta^{\textit{ab}}(R)$,
  the free commutative monoid generated by the $\partial_i$'s. It is isomorphic to $\N^{d}$. Every element $\theta$ of $\Theta^{\textit{ab}}(R)$ can be written uniquely
\[ \theta=\prod_{i=1}^d {\partial_i}^{e_i(\theta)}. \] 
They act in a natural way on $R$. The integer $e(\theta):=\sum_i e_i(\theta)$ is called the \emph{order of $\theta$}. It satisfies $e(\theta_1 \theta_2)=e(\theta_1)e(\theta_2)$. 
For any $\theta, \theta'\in\Theta^{\textit{ab}}(R)$, if $\theta'$ divides $\theta$ (\ie if $e_i(\theta')\leq e_i(\theta)$ for all $i$), Kolchin defines
\[ \binom{\theta}{\theta'}:= \prod_{i=1}^d \binom{e_i(\theta)}{e_i(\theta)'}. \]
\smallskip

\begin{proposition}[\protect{\Cite[p. 60]{DiffAlgGroupsKolch}}]\label{Prop_Pascal_Kolchin}
Let $R$ be a partial \D-ring. Then, for every $f,g\in R$ and every $\theta\in \Theta^{\textit{ab}}(R)$, one has
\[ 
\theta (fg)=\sum_{\begin{subarray}{c}
 (\theta_1, \theta_2)  \\
s.t.\  \theta_1 \theta_2=\theta
\end{subarray} } \binom{\theta}{\theta_1}\cdot \theta_1(f) \cdot\theta_2 (g).
 \]
\end{proposition}

\pproof{
First, let us remark that if one sets $\binom{\theta}{\theta'}:=0$ whenever $\theta'$ does not divide $\theta$, then  these generalized binomial numbers  satisfy a Pascal-like identity:
\[ 
\fa i_0\in \lbrace 1, \ldots, d\rbrace, \quad \fa \theta, \theta'\in \Theta^{\textit{ab}}(R), \quad
\binom{\partial_{i_0}\theta}{\partial_{i_0}\theta'}=\binom{\theta}{\theta'}+\binom{\theta}{\partial_{i_0}\theta'}.
 \]
Now, we proceed by induction on $e(\theta)$. When $e(\theta)=1$, this formula reduces to the Leibniz rule. Let us assume it is true if $e(\theta)\leq n$, and let $\theta$ be a differential operator of order $n$. Let $\partial$ be any of $\partial_i$'s. One has
\begin{align*}
\partial \theta (fg)&=\sum_{\theta_1 \theta_2 = \theta}  \binom{\theta}{\theta_1}\cdot \partial\theta_1(f) \cdot\theta_2 (g)+ \binom{\theta}{\theta_1}\cdot \theta_1(f) \cdot\partial\theta_2 (g)\\
&=\sum_{
\begin{subarray}{c}
\widetilde{\theta}_1 \widetilde{\theta}_2 =\partial \theta \\
\partial \mid \widetilde{\theta}_1
\end{subarray}
}
\binom{\theta}{\widetilde{\theta}_1/\partial} \cdot \widetilde{\theta}_1(f)\cdot \widetilde{\theta}_2(g)
+
\sum_{
\begin{subarray}{c}
\widetilde{\theta}_1 \widetilde{\theta}_2 =\partial \theta \\
\partial \mid \widetilde{\theta}_2
\end{subarray}
}
\binom{\theta}{\widetilde{\theta}_2/\partial} \cdot \widetilde{\theta}_1(f)\cdot \widetilde{\theta}_2(g) \\
&=\sum_{
\begin{subarray}{c}
\widetilde{\theta}_1 \widetilde{\theta}_2 =\partial \theta \\
\partial \mid \widetilde{\theta}_1
\end{subarray}
}
\binom{\theta}{\widetilde{\theta}_1/\partial} \cdot \widetilde{\theta}_1(f)\cdot \widetilde{\theta}_2(g)
+
\sum_{
\begin{subarray}{c}
\widetilde{\theta}_1 \widetilde{\theta}_2 =\partial \theta \\
\partial \mid \widetilde{\theta}_2
\end{subarray}
}
\binom{\theta}{\widetilde{\theta}_1} \cdot \widetilde{\theta}_1(f)\cdot \widetilde{\theta}_2(g)\
\\
&=
\sum_{
\begin{subarray}{c}
\widetilde{\theta}_1 \widetilde{\theta}_2 =\partial \theta \\
\partial \mid \widetilde{\theta}_1 \text{ and } \partial \mid \widetilde{\theta}_2
\end{subarray}
}
\left ( \binom{\theta}{\widetilde{\theta}_1/\partial}+\binom{\theta}{\widetilde{\theta}_1}\right) \cdot \widetilde{\theta}_1(f)\cdot \widetilde{\theta}_2(g) \\*
& \qquad \qquad +\sum_{
\begin{subarray}{c}
\widetilde{\theta}_1 \widetilde{\theta}_2 =\partial \theta \\
\partial \mid \widetilde{\theta}_1 \text{ and } \partial \nmid \widetilde{\theta}_2
\end{subarray}
}
 \binom{\theta}{\widetilde{\theta}_1/\partial} \cdot \widetilde{\theta}_1(f)\cdot \widetilde{\theta}_2(g) 
\\*
& \qquad \qquad  \qquad \qquad+
 \sum_{
\begin{subarray}{c}
\widetilde{\theta}_1 \widetilde{\theta}_2 =\partial \theta \\
\partial \mid \widetilde{\theta}_2 \text{ and } \partial \nmid \widetilde{\theta}_1
\end{subarray}
}
 \binom{\theta}{\widetilde{\theta}_1} \cdot \widetilde{\theta}_1(f)\cdot \widetilde{\theta}_2(g) 
\end{align*}
But, when  $\partial\nmid \widetilde{\theta}_2$ and $\widetilde{\theta}_1 \widetilde{\theta}_2 =\partial \theta$, then $\widetilde{\theta}_1 \nmid \theta$. Thus one has, in this case,
\[ 
\binom{\theta}{\widetilde{\theta}_1/\partial}=\binom{\theta}{\widetilde{\theta}_1/\partial}+\binom{\theta}{\widetilde{\theta}_1}
 \]
 and a similar formula when $\partial\nmid \widetilde{\theta}_1$ and $\widetilde{\theta}_1 \widetilde{\theta}_2 =\partial \theta$. Thus, 
\begin{align*}
\partial \theta (fg)&=
\sum_{
\widetilde{\theta}_1 \widetilde{\theta}_2 =\partial \theta 
}
\left ( \binom{\theta}{\widetilde{\theta}_1/\partial}+\binom{\theta}{\widetilde{\theta}_1}\right) \cdot \widetilde{\theta}_1(f)\cdot \widetilde{\theta}_2(g) 
\\&=
\sum_{
\widetilde{\theta}_1 \widetilde{\theta}_2 =\partial \theta 
}
 \binom{\partial\theta}{\widetilde{\theta}_1}\cdot \widetilde{\theta}_1(f)\cdot \widetilde{\theta}_2(g) 
\end{align*}
by Pascal's identity. }

\begin{fact}
The monoid $\Theta^{\textit{ab}}(R)$ has  well-orders.
\end{fact} 

\pproof{As already said, $\Theta^{\textit{ab}}(R)$ is isomorphic to $\N^d$. So, the lexicographical order can be transported to $\Theta^{\textit{ab}}(R)$. It is a compatible (for the monoid structure) well-order.
\smallskip

Here is another compatible well-order for $\Theta^{\textit{ab}}(R)$. By considering the map (see \Cite{KeigherDiffRingFromQuasPrimes} or \Cite{RittDiffAlg})
\[ \phi 
:
\fonctionb{\Theta^{\textit{ab}}(R)}{\N\times\N^d}{\theta}{e(\theta), (e_i(\theta))_i}
 \]
 one endows $\Theta^{\textit{ab}}(R)$ with an order, defined by $\theta < \theta'$ iff $\phi(\theta)\prec_{\text{lex}} \phi(\theta')$. The verifications are left to the reader. Note than Robbiano gives in \Cite{RobbianoTermOrdering} a classifiation of all total orders on $\N^{d}$ compatible with the monoid structure. }
 
\smallskip

\begin{proposition}\label{P_diese_premier_bis}
Let $R$ be partial $\Q$-\D-algebra and let $\mathfrak{p}$ be a prime ideal of $R$. Then,
$ \mathfrak{p}_\#$
is prime.
\end{proposition}

\pproof{This proof is after W. Keigher (see Proposition 1.5 of \Cite{KeigherDiffRingFromQuasPrimes}), but generalized to any compatible well-order on $\Theta^{\textit{ab}}(R)$. Let $\mathfrak{p}$ be a prime ideal of $R$ and let $x,y\in R\setminus \mathfrak{p}$.  Let $<$ be a compatible well-order on  $\Theta^{\textit{ab}}(R)$ so that we can define
\[
\theta_*:= \min \left\lbrace \theta \in \Theta(R)\tq \theta(*)\notin \mathfrak{p} \right\rbrace 
 \]
 for $*=x,y$. Now, by Proposition \Ref{Prop_Pascal_Kolchin},
 \[ 
(\theta_x \theta_y) (xy) = \binom{\theta_x \theta_y}{\theta_x} \theta_x(x)\theta_y(y)+
\sum_{\begin{subarray}{c}
\theta\theta'=\theta_x \theta_y \\
\theta\neq \theta_x, \theta'\neq \theta_y
\end{subarray}}
 \binom{\theta_x \theta_y}{\theta} \theta(x)\theta'(y).
  \]
In the indexed sum, one cannot have $\theta>\theta_x$ and $\theta'>\theta_y$, since $>$ is compatible with product. So, for instance $\theta\leq\theta_x$, but since  $\theta\neq\theta_x$, one has $\theta<\theta_x$. So, $\theta(x)\in \mathfrak{p}$ --- and the the same is true for the big sum. Now, since $\mathfrak{p}$ is prime and since  $ \binom{\theta_x \theta_y}{\theta_x}$ is invertible in $R$, one has $\binom{\theta_x \theta_y}{\theta_x} \theta_x(x)\theta_y(y)\notin \mathfrak{p}$, and so $(\theta_x \theta_y) (xy)\notin \mathfrak{p}$. Hence, $xy\notin \mathfrak{p}_\#$ and thus $\mathfrak{p}_\#$ is prime.
}

\bigskip

\subsection{Primality of \D-ideals in partial \D-rings}\label{Prim_D_ideals}

In this section, we generalize to \D-rings some classical results on primality of \D-ideals. As a byproduct, we are able to prove Proposition \Ref{P_diese_premier_bis} without any assumption of  commutativity or finiteness. See  \Cite{kaplansky} for a nice exposition of some of these results in the case of ordinary \D-rings. Lemmas \Ref{lemma_easy}, \Ref{lemma_nilradical} and Proposition \Ref{radical_diff} are immediate when one knows the similar results for ordinary \D-rings. Only Proposition \Ref{primality_ideals} requires a little bit of work.

\begin{lemma}\label{lemma_easy}
Let $I$ be a \D-ideal. Let $\theta\in \Theta(R)$. Then, $$\fa x,y\in R, \qquad xy\in I\ \impl \ x^{e(\theta)+1}\cdot \theta(y)\in I$$
\end{lemma}

\pproof{By induction on $e(\theta)$. See also \Cite[p. 62]{DiffAlgGroupsKolch} in the partial case. }
\begin{lemma}\label{lemma_nilradical}
Let $R$ be a $\Q$-algebra. Let $x\in R$ such that $x^n=0$. Then, for all $\partial\in \mathrm{Der}(R)$, and all $\ell\in \{1, \ldots, n\}$, one has
\[ 
(\partial x)^{2 \ell-1}\cdot x^{n-\ell}=0.
 \]
 In particular, $(\partial x)^{2n-1}=0$
\end{lemma}

\pproof{Straightforward by induction on $\ell$. The hypothesis that $R$ is a $\Q$-algebra in highly used there. See \Cite{kaplansky}, Lemma 1.7, page 12. }

\smallskip

\begin{proposition}\label{radical_diff} Let $R$ be a $\Q$-\D-algebra and $I$ a \D-ideal of $R$. Then,
\begin{itemize}\setlength{\itemsep}{0mm}
\item[(i)] The nilradical $\Nil(R)$ of $R$ is a \D-ideal.

\item[(ii)] More precisely, if $x\in R$ and $x^{n}=0$, then for all $\theta\in \Theta(R)$, one has
$\theta(x)^{2^{e(\theta)}(n-1)+1} =0$
\item[(iii)] The radical $\sqrt{I}$ is a \D-ideal. 
\end{itemize}
\end{proposition}

\pproof{One gets \emph{(iii)} by applying  \emph{(i)} for $R/I$. Now,  \emph{(i)} is a consequence of \emph{(ii)}, and \emph{(ii)} is a  consequence of Lemma \Ref{lemma_nilradical}.
}

\smallskip

\begin{proposition}\label{primality_ideals}
Let $R$ be a $\Q$-\D-algebra. Let $I$ be a \D-ideal of $R$ and let $\mathfrak{p}$ be a prime ideal of $R$ such that $I \subset \mathfrak{p}$.
\begin{itemize}\setlength{\itemsep}{0mm}
\item[(i)] Every maximal proper  \D-ideal of $R$ containing $I$ is prime.
\item[(ii)] Every minimal prime ideal of $R$ is a \D-ideal.
\item[(iii)] Every maximal \D-ideal $J$ such that $I\subset J\subset \mathfrak{p}$ is prime. 
\item[(iv)] The \D-ideal $\mathfrak{p}_\#$ is prime.
\end{itemize}
\end{proposition}

\bigskip

Let us recall that for any ideal $I$, $I_\#$ is defined by 
\[ 
I_\#:= \left\lbrace x\in R \tq \fa \theta \in \Theta(R), \quad\theta(x)\in I\right\rbrace .
 \]
Before proving this proposition, let us give counterexamples to Proposition \Ref{radical_diff} and \Ref{primality_ideals} when  $R$ is not defined over $\Q$. First, take a field $k$ of characteristic $p>0$, and consider the ordinary \D-ring $k[x]$ with $x'=1$. Denote $I:=(x^{p})$. Then,  $I$ is a \D-ideal, and is maximal among  proper \D-ideals, but is not prime. Furthermore, $\sqrt{I}=(x)$ is not a \D-ideal (but is prime). This gives counterexamples to Proposition \Ref{radical_diff} and  \emph{(i)}, \emph{(iii)}  of Proposition \Ref{primality_ideals}. Second, since $I$ is a \D-ideal, one can consider the \D-ring $R:= k[x]/I$. In $R$, the ideal $(x)$ is prime and minimal --- but is not a \D-ideal. Furthermore, one has $(x)_{\#}=(0)$, which is not prime. This gives counterexamples to \emph{(ii)} and \emph{(iv)} of Proposition \Ref{primality_ideals}.
\bigskip

\pproof{
First, \emph{(iv)} is a consequence of \emph{(iii)}, with $I=0$. Indeed, one can easily check that $\mathfrak{p}_\# $ is the largest  prime \D-ideal contained in $\mathfrak{p}$. Now, let us prove \emph{(i)} and \emph{(iii)}. They both will be a particular case of the following result. 
\begin{center}
\begin{itshape}
 Let $S$ be a multiplicative subset of $R$ and let $K$ be a \D-ideal of $R$ such that $S \cap K=\varnothing$. Let us consider $L_0$ a \D-ideal maximal amongst those that verify $K \subset L$ and $L \cap S=\varnothing$. Such an $L_0$ exists by Zorn's lemma. 
Then, $L_0$ is prime. 
 \end{itshape}
\end{center}First, $\sqrt{L_0}=L_0$. Indeed, by Proposition~\Ref{radical_diff}.\emph{(iii)}, $\sqrt{L_0}$ is a \D-ideal and one checks that  $\sqrt{L_0}\cap S=\varnothing$ and $K\subset \sqrt{L_0}$. By maximality, $\sqrt{L_0}=L_0$. Now, let us prove that for all $s\in S$, the ideal
$(L_0 : s^\infty)$ (see paragraph \RefPar{colon_ideals})
equals $L_0$. This ideal $(L_0 : s^\infty)$, by Lemma \Ref{lemma_colon_ideal},  contains $L_0$ and is a \D-ideal. Furthermore, one has $S\cap (L_0:s^{\infty})=\varnothing$.  Indeed if $t\in S \cap (L_0:s^{\infty})$, that means there exists $m>0$ such that $ts^m\in L_0$, which is absurd. So, by maximality, $L_0=(L_0 : s^\infty)$.

\smallskip

 Now, let $x,y\in R$ with $xy \in L_0$ and $y\notin L_0$. As for  $(L_0 : s^\infty)$, the ideal $(L_0 : y^\infty)$ contains $L_0$ and is a \D-ideal. Let us prove that it does not intersect $S$.
 Suppose the contrary and let $s\in S$ such that there exists $n>1$ such that $sy^{n}\in L_0$. Since $L_0=\sqrt{L_0}$, one also has $sy\in L_0$. This means $y\in (L_0 : s^{\infty})$. As this latter equals $L_0$ it is absurd. So, by maximality, one has $L_0 = (L_0 : y^{ \infty})$. Since, $x\in (L_0 : y^{ \infty})$ one has $x\in L_0$. Hence, $L_0$ is prime. 
 
\smallskip

Let us come back to the proposition. One obtains \emph{(i)} by applying this result with $S=\{1\}$ and $K=\{0\}$.  One obtains \emph{(iii)} by applying this result with $S=R\setminus \mathfrak{p} $ and $K=I$.
\smallskip

Finally, let us prove \emph{(ii)}. Let $\mathfrak{p}$ be a minimal prime ideal (the existence of such ideal is guaranteed by Zorn's lemma). Then, by \emph{(iv)}, $\mathfrak{p}_\#$ is also prime. So, by minimaliy,   $\mathfrak{p}_\#=\mathfrak{p}$, and $\mathfrak{p}$ is a \D-ideal. }

\bigskip

\subsection{Simple \D-rings}\label{section_simple_D_rings}

The aim of this paragraph is to provide a reference for properties of simple \D-rings in the general case. Some of the results presented are known and proved for ordinary \D-rings --- we generalize it to \D-rings. We also prove some new facts (see point \emph{(i)} of Proposition \Ref{Properties_simple}), and an interesting lemma (see Lemma \Ref{Super_lemme}). The point \emph{(iv)} is proved in \Cite{singvdp} for simple ordinary \D-rings whose field of constants is algebraically closed of characteristic zero. We prove it for arbitrary simple \D-rings, and replace in the proof the use of Chevalley's theorem on constructible sets by the application of Noether normalization lemma. First, recall:

\begin{definition}
A \D-ring $R$ is said to be \emph{simple} if the only \D-ideals of $R$ are $(0)$ and $R$.
\end{definition}

\begin{proposition}\label{Properties_simple}
Let $R$ be a simple \D-ring. Then,
\begin{itemize}\setlength{\itemsep}{0mm}
\item[(i)] \begin{itemize}\setlength{\itemsep}{0mm}
\item[a)] $R$ is irreducible.
\item[b)] $R$ is connected.
\end{itemize} 
\item[(ii)] The constant ring $C(R)$ of $R$ is a field.
\item[(iii)] If $R$ is a $\Q$-algebra  or if $R$ is reduced, then $R$ has no zero divisors.
\item[(iv)] Let $S$ be multiplicative subset of $R$ not containing zero. If we denote by $i_S : R \longto S^{-1}R$ the localization map, then:
\begin{itemize}
\item[a)] the morphism of rings $C(i_S) : C(R)\longto C(S^{-1}R)$ is an isomorphism;
\item[b)] the \D-ring $S^{-1}R$ is simple.
\end{itemize}
\item[(v)] Suppose furthermore that $R$ is a $k$-\D-algebra, finitely generated and without zero divisors, for a \D-field $k$.  Then, the extension $C(R)/C(k)$ is algebraic.
\item[(vi)] Let $L$ be a constant field extension of $C(R)$. Then, $R\otimes_{C(R)} L$ is simple, and $C(R\otimes_{C(R)} L)=L$.
\end{itemize}
\end{proposition}

\smallskip

The point \emph{(iii)} is false in general. Indeed, the  ordinary \D-ring $k[x]/x^p$, for any field $k$ of characteristic $p>0$, with the derivation $x'=1$ is simple.

\smallskip

To prove \emph{(v)}, we need two lemmas:

\begin{lemma}\label{lemma_alg_cst}
Let $k$ be a \D-field and $R$ a $k$-\D-algebra. Let $f\in C(R)$ such that $f$ is algebraic over $k$. Then, $f$ is algebraic over $C(k)$.
\end{lemma}

\pproof{
Consider the minimal polynomial $P\in k[X]$ of $f$, differentiate $P(f)=0$ and conclude.
}

\begin{lemma}\label{Lemma_image_element_inv}
 Let $k$ be a field and $k\to \overline{k}$ an algebraic closure of $k$. Let $R$ be a finitely generated $k$-algebra, without zero divisors. Let $f\in R^\times$. Let us denote
 \[ 
 \Im_{\overline{k}}(f):=\left \{ 
 \varphi(f) \tq \varphi\in \Hom_{\Alg{k}}(R, \overline{k})
 \right \}.
 \]
 Then: --- if $f$ is algebraic over $k$, $\Im_{\overline{k}}(f)$ is finite ; \\
\phantom{Then: }--- otherwise, $\Im_{\overline{k}}(f)=\overline{k}\setminus \{0\}$ 
\end{lemma}
\smallskip

\pproof{
In \Cite{singvdp}, a slighly different  form of this lemma is given, deduced from Chevalley's theorem on constructible sets. We give here a refinement and a different proof. If $f$ is algebraic and if $P\in k[X]$  is such that $P(f)=0$, then for all $\varphi\in\Hom_{\Alg{k}}(R, \overline{k})$, $P(\varphi(f))=0$, and so $ \Im_{\overline{k}}(f)$ is included of the set of roots of $P$ in $\overline{k}$.

\smallskip
Now, let us set $K=k(f)$. $R$ is a finitely generated $K$-algebra without zero divisors.  By Noether's normalization lemma, one can thus find $x_1, \ldots, x_n$ such that the $x_i$ are algebraically independent over $K$ and $R$ is integral over $K[x_1, \ldots, x_n]$. Let $a\in \overline{k}$, $a\neq 0$. If $f$ is transcendent over $k$, one can find a morphism $\varphi : K \longto \overline{k}$ such that $\varphi(f)=a$. One can extend $\varphi$ to $K[x_1, \ldots, x_n]$ for instance by setting $\varphi(x_i)=0$.  Now, since $R/ K[x_1, \ldots, x_n]$ is integral, by Proposition 3.1 of chap VII of \Cite{algebra}, one knows that there exist an extension of $\varphi$ to $R$. Hence, $a\in \Im_{\overline{k}}(f)$.
}

\bigskip

To prove item \emph{a)} of \emph{(ii)}, we need the following lemma, which is also interesting by itself.

\begin{lemma}\label{Super_lemme}
Let $R$ be a ring. Let $I$ and $J$ be two ideals of $R$ such that $I \cap J=0$. Let $x\in I$ and $y\in J$. Then, 
\[\fa \theta, \theta'\in \widehat{\Theta}(R), \quad\theta(x)\cdot \theta'(y)=0.\]
\end{lemma}

\smallskip

\pproof{
We prove it by induction on $e(\theta)+e(\theta')$. If $e(\theta)+e(\theta')=0$, this reduces to prove that $xy=0$. But $xy \in I \cap J$. Now, assume that the result holds for all $\theta, \theta'\in\Theta$ such that $e(\theta)+e(\theta')=N$. Let $\theta_0, \theta_1\in\Theta$ such that $e(\theta_0)+e(\theta_1)=N+1$. More precisely, we set $e(\theta_0)=i$ and $e(\theta_1)=N+1-i$.
\smallskip

 In what follows, we will use this notation: given a differential operator $\theta:=\delta_1 \delta_2 \cdots \delta_p$, we will denote
\[ 
[\, \theta\, ]_k^{\ell}:= \delta_k  \delta_{k+1} \cdots \delta_\ell
 \]
if $k\leq \ell$ and $[\theta]_k^{\ell}=\Id$ if $k >\ell$. Thus, the order of $[\, \theta\, ]_k^{\ell}$ is 
\[ e \left ( [\, \theta\, ]_k^{\ell} \right )=\max (0, \ell-k+1).\]
Furthermore, we will write
$\overline{\theta}:=\delta_p\delta_{p-1}\cdots \delta_1$.

\smallskip

 By assumption, one has for every $j\in \{2, \ldots, i+1\}$
\begin{equation}\label{hyp_rec_formelle}
 [\, \theta_0 \, ]_j^{i} (x)  \cdot  ( [\, \overline{\theta_0\, ]_1^{j-2}}   \ \theta_1) (y)=0
\end{equation}  
  as one easily checks that the sum of the orders of the involved differential operators is equal to $N$ in all the cases. Applying the $(j-1)$-th derivation in the writing of $\theta_0$, that is, applying the derivation $ [\, \theta_0 \, ]_{j-1}^{j-1}$ to (\ref{hyp_rec_formelle}), one gets
\begin{equation}\label{formelle_2}
 [\, \theta_0 \, ]_{j-1}^{i} (x)  \cdot  ( [\, \overline{\theta_0\, ]_1^{j-2}}   \ \theta_1) (y) +  [\, \theta_0 \, ]_{j}^{i} (x)  \cdot  ( [\, \overline{\theta_0\, ]_1^{j-1}}   \ \theta_1) (y)=0.
\end{equation} 
Denoting
\[ 
\alpha_{j}:= [\, \theta_0 \, ]_{j-1}^{i} (x)  \cdot  ( [\, \overline{\theta_0\, ]_1^{j-2}}   \ \theta_1) (y) 
 \]
 one can write (\ref{formelle_2}) as $\alpha_j+\alpha_{j+1}=0.$ So,
 \[ 
\sum_{j=2}^{i+1} (-1)^{j} (\alpha_j+\alpha_{j+1}) =\alpha_2+(-1)^{i+1}\alpha_{i+2}=0.
  \]
But, 
\[ 
\alpha_2=\theta_0(x)\cdot  \theta_1(y)\quad \aand \quad \alpha_{i+2}=x \cdot (\overline{\theta_0} \ \theta_1)(y).
 \]
 Hence, 
 \begin{align}
 \theta_0(x)\cdot  \theta_1(y)+ (-1)^{i+1}\cdot  x \cdot (\overline{\theta_0} \ \theta_1)(y)&=0\label{equation_premiere} \\
  \theta_1(y) \cdot  \theta_0(x) + (-1)^{(N+1-i)+1}\cdot  y \cdot (\overline{\theta_1} \ \theta_0)(x)&=0, \label{equation_duale}
\end{align}
the equation (\ref{equation_duale}) being obtained by
 interchanging $x\leftrightarrow y $ and $\theta_0\leftrightarrow \theta_1$. Thus, one gets
 \begin{equation}
(-1)^{N}y \cdot (\overline{\theta_1} \ \theta_0)(x) - x \cdot (\overline{\theta_0} \ \theta_1)(y)=0
\end{equation}
 But, this implies that $x \cdot (\overline{\theta_0} \ \theta_1)(y)\in I \cap J$ and thus is zero. Thus, by (\ref{equation_premiere}), one has 
 $ \theta_0(x)\cdot  \theta_1(y)=0$, what we wanted. 
}

\bigskip

\pproofbis{of Proposition \Ref{Properties_simple}}{
For \emph{(ii)}, \emph{(iv)} and \emph{(vi)}, see respectively Theorem 2, Theorem 6 and Theorem 7 of \Cite{DiffRings}. 
\smallskip

The proof of item \emph{b)} of \emph{(i)} is  easy. It is a direct consequence of item \emph{a)}, but let us give a direct proof. Let $e\in R$ be a nonzero idempotent. Let $'$ be any derivation of $R$. Since $e^{2}=e$, one has $2ee'=e'$ and so $2ee'=ee'$ and so $ee'=0$ and so $e'=0$. Thus, $(e)$ is a \D-ideal. Hence, $e$ is invertible and so $e=1$.
\smallskip

Let us prove now item \emph{a)} of \emph{(i)}. Let $I$ and $J$ be two nonzero ideals. We want to prove that $I \cap J\neq 0$. Assume $I \cap J=0$. Let $\theta\in \Theta(R)$. We denote by $\theta(J):=\{\theta(y), y\in J\}$. By Lemma \Ref{Super_lemme}, one easily sees that $\theta(J)\cdot I=0$ and thus $\theta(J)\subset \mathrm{Ann}(I)$, the annihilator of $I$. But, since $J$ is nonzero and $R$ simple, the ideal generated by all the $\theta(J)$ equal $R$. Hence, $\mathrm{Ann}(I)=R$ and $I=0$, which is absurd. 

\smallskip

Let us now prove \emph{(iii)} and assume $R$ reduced. Let $x,y\in R$ such that $xy=0$, and assume $x\neq 0$. By Lemma  \Ref{lemma_colon_ideal}, the ideal $(0:y^\infty)$ is a \D-ideal. But $x\in(0:y^\infty)$ and so $(0:y^\infty)=R$ and so there exists $n\geq 1$ such that $y^n=0$ and so $y=0$. If $R$ is a $\Q$-algebra, then $\Nil(R)$ is a \D-ideal and so $R$ is reduced. 

\smallskip

Let us now prove \emph{(v)}. Let $f\in C(R)$, $f\neq 0$. By \emph{(ii)}, we know that $f$ is invertible. Let $k\to \overline{k}$ be an algebraic closure of $k$. Let $\varphi\in \Hom_{\Alg{k}}(R, \overline{k})$. But, since $f-1$ lies also in $C(R)$, it is also invertible and so, $\varphi(f-1)$ cannot be zero, and so, $\varphi(f)\neq 1$. Thus,  by Lemma \Ref{Lemma_image_element_inv}, $f$ is algebraic over $k$. But, since $f$ is constant, one knows by Lemma \Ref{lemma_alg_cst} that $f$ is algebraic over $C(k)$.
}

\bigskip

One also has the following characterization of simple \D-rings:

\begin{fact}\label{char_simple}
Let $R$ be a $\Q$-\D-algebra whose only  prime \D-ideal of $R$ is $(0)$. Then, $R$ is simple.
\end{fact}

This fact is false when $R$ is not defined over $\Q$. Indeed, consider $R:=k[x]$ where $k$ is a field of characteristic $p>0$, with $x'=1$.  The only prime \D-ideal of $R$ is $(0)$ but $(x^p)$ is a \D-ideal of $R$.

\smallskip

\pproof{
Let $I$ be a \D-ideal. Then, by Proposition \Ref{primality_ideals}.\emph{(i)}, $I$ is included in a prime \D-ideal. Hence, $I=(0)$.
}

\bigskip

To end up, let us remark  the following consequence of Lemma \Ref{Super_lemme}. We need to introduce two notations. If $R$ is a ring and if $I$ is an ideal of $R$, we denote by
\[
[I ]:=\sum_{\theta \in \widehat{\Theta}(R)} \theta (I)
\] 
the ideal generated by the images of $I$ under all the possible differential operator of $R$. It is the smallest ideal of $R$ containing $I$ and stable under all derivations of $R$. We also set  
\[
\mathrm{Nil}_n(R):=\{ x\in R\tq x^{n}=0\},
\]
for $n\in \Z_{\geq 0}$.

\smallskip

\begin{corollary}\label{Coro_Ttes_deriv}
Let $R$ be a ring and let $I,J$ ideals such that $I\cap J=0$. Then, $[I]\cap [J]\subset \mathrm{Nil}_2(R)$.
\end{corollary}

\pproof{
It is very easy. Assume that one $x\in [I]\cap [J]$ and write 
\[
x=\sum_{i\in I} \theta_i(a_i)=\sum_{j\in J} \theta_j(b_j)
\]
with $\theta_i, \theta_j\in \widehat{\Theta}(R)$ and $a_i\in I$, $b_j\in J$. Then,
\[
x^{2}=\sum_{(i,j)\in I\times J} \theta_i(a_i)\cdot \theta_j(b_j)
\]
which is equal to $0$ by Lemma \Ref{Super_lemme}.
}

\bigskip 

Here is a counterexample to $[I]\cap[J]=0$ when $I\cap J=0$. Take $R=k[x,y]/(x^{2}, xy, y^{2})$ with $k$ a field of characteristic $2$ and consider the derivation $\partial$ defined by $\partial x=y$ and $\partial y=0$. Set $I:=(x)$ and $J:=(y)$. Then, one has $I\cap J=0$ but $\langle I \rangle=(x,y)$. Hence, $[I]\cap[J]=0$  can not hold. 

\smallskip

The geometric interpretation of  Corollary \Ref{Coro_Ttes_deriv} is the following. The proof is left to the reader.

\smallskip

\begin{corollary}\label{Coro_geo_irr}
Let $X$ be a scheme and let $F_1, F_2$ be two subsets of $X$ such that $F_1\cup F_2=X$. Define, for any closed set $F$,
\[
[F]:= \{ x\in F \tq \fa \chv{V}\in \mathscr{T}_X(X), \ \mathrm{Traj}_{\chv{V}}(x)\in F\}.
\]
Then, $[F_1]\cup[F_2]=X$.
\end{corollary}
\newpage

\def\cprime{$'$}
\providecommand{\bysame}{\leavevmode ---\ }
\providecommand{\og}{``}
\providecommand{\fg}{''}
\providecommand{\smfandname}{\&}
\providecommand{\smfedsname}{\'eds.}
\providecommand{\smfedname}{\'ed.}
\providecommand{\smfmastersthesisname}{M\'emoire}
\providecommand{\smfphdthesisname}{Th\`ese}

\end{document}